\documentclass[11pt]{amsart}
\usepackage{mathcomp,amscd,xypic}
\usepackage{amssymb, amscd, amsthm}
\usepackage{times}
 \usepackage{xypic}
 \usepackage{upquote}
\usepackage{hyperref}
\usepackage{verbatim}
\usepackage{color} 
\usepackage{alltt}
\DeclareFontEncoding{OT2}{}{} %

 \newcommand{\OLrho}{\mathrm{U}_g}
 
\newcommand{\GL}{\mathrm{GL}}
\newcommand{\et}{\mathrm{et}}
\newcommand{\vM}{\varpi_{\mathrm{Merel}}}
\newcommand{\mot}{\mathrm{mot}}
\newcommand{\UnitL}{\mathrm{U}_L}

\newcommand{\End}{\mathrm{End}}
\newcommand{\ind}{\mathrm{ind}}
\newcommand{\Frob}{\mathrm{Frob}}
\newcommand{\res}{\mathrm{res}}

\newcommand{\Shim}{\mathfrak{S}}

\newcommand{\Int}{\mathrm{int}}
\newcommand{\Gproj}{G^{\mathrm{proj}}}
\newcommand{\Zar}{\mathrm{Zar}}

\renewcommand{\dim}{\mathrm{dim}}

\newcommand{\PP}{\mathbf{P}}
\newcommand{\QQ}{\mathbf{Q}}
\newcommand{\rhobar}{\overline{\rho}}

\newcommand{\Ext}{\mathrm{Ext}}

\newcommand{\Spec}{\mathrm{Spec}}
\newcommand{\ur}{\mathrm{ur}}

\newcommand{\Ad}{\mathrm{Ad}}
\newcommand{\Mg}{\mathrm{M}_g}

\newcommand{\fl}{\mathrm{fl}}

\newcommand{\Sh}{\mathbf{S}}

\newcommand{\cFd}{\mathcal{F}^{\bullet}}
\newcommand{\cHd}{\mathcal{H}^{\bullet}}

\newcommand{\ra}{\rightarrow}
\newcommand{\isoarrow}{{~\overset\sim\longrightarrow}}

\usepackage[OT2,T1]{fontenc}
\DeclareSymbolFont{cyrletters}{OT2}{wncyr}{m}{n}

\newcommand{\C}{{\mathbb C}}
\newcommand{\Z}{{\mathbb Z}}

\newcommand{\F}{{\mathbf F}}

\newcommand{\Q}{{\mathbf Q}}
 
\newcommand{\Gal}{\mbox{Gal}}

\newcommand{\Hom}{\mathrm{Hom}}

  \newcommand{\cA}{\mathcal{A}}

\newcommand{\CL}{\omega}

\newcommand{\disc}{\mbox{disc}}

\numberwithin{equation}{section}
\numberwithin{table}{section}
\numberwithin{figure}{section}

\newtheorem{Proposition}[subsection]{Proposition} 
\newtheorem{Lemma}[subsection]{Lemma}

\newtheorem{Conjecture}{Conjecture}[section]
\newtheorem{lemma}[subsection]{Lemma}
\newtheorem{rmk}{Remark}[section]
\allowdisplaybreaks
 \swapnumbers

\begin{document}

\author{Michael Harris}
\thanks{M.H.'s research received funding from the European Research Council under the European Community's Seventh Framework Programme (FP7/2007-2013) / ERC Grant agreement no. 290766 (AAMOT).  M.H. was partially supported by NSF Grant DMS-1404769}
\author{Akshay Venkatesh}
\thanks{A. V. research was partially supported by the NSF and by the Packard foundation}

\title{Derived Hecke algebra for weight one forms}

\begin{abstract}
 We study the action of the derived Hecke algebra on the space of weight one forms. By analogy
with the topological case, we formulate a conjecture relating this to a  certain Stark unit.

We verify the truth of the conjecture numerically,  for the weight one forms of level $23$ and $31$, and  many derived Hecke operators at primes less than $200$.   Our computation
depends in an essential way on Merel's evaluation of the pairing between the Shimura and cuspidal subgroups of $J_0(q)$. 
\end{abstract}

\maketitle

\tableofcontents
  
\section{Introduction}

Let $G$ be an algebraic group over $\QQ$. In \cite{V3} the second-named author
studied the action of a derived version of the Hecke algebra on the singular cohomology  of the locally symmetric space   attached to $G$.
One expects that this action transports Hecke eigenclasses between cohomological degrees and 
moreover (see again \cite{V3}) is related to a ``hidden'' action of a motivic cohomology group.

 It is also possible for  a Hecke eigensystem  on {\em coherent} cohomology to occur in multiple degrees.
The simplest situation is 
weight one forms for the modular curve.   We study this 
case,  explicating the action of the derived Hecke algebra and 
formulating a conjectural relationship with motivic cohomology.

The motivic cohomology is particularly concrete: a weight one eigenform $f$ is attached to a $2$-dimensional Artin representation $\rho_f$ of $\mathrm{Gal}(\bar{\Q}/\Q)$ and the motivic  cohomology group
in question  is   generated by a certain unit in the splitting field of the adjoint of  $\rho_f$.  
This makes the conjecture particularly amenable to numerical testing. We  carry this out for the forms
of conductor $23$ and $31$, and the first few derived Hecke operators; the numerics all support the conjecture. 

We note that this story is  also related to the Taylor--Wiles method for coherent cohomology  and its obstructed version due to Calegari--Geraghty (see   \cite{H}, and the  detailed discussion of weight one forms in \cite{CG}). 
This is used implicitly in the discussion in \S \ref{GDS} of the current paper.

 Now we outline the contents of the paper. After giving notation in \S \ref{notn}, we  describe the derived Hecke algebra and the main conjecture  in \S \ref{DHmain}.  Although the discussion to this point is self-contained,
 we postpone the comparison with the results of \cite{V3} until the final section. We translate the Conjecture to an explicitly computable form
 in \S \ref{sec:explicit}; see in particular Proposition \ref{prop:explicit}. Finally, in \S \ref{sec:cubic} we make the conjecture even more explicit
 in the case of a form associated to a cubic field $K$, and check it numerically in the case of $K$ with discriminant $-23$ and $-31$. 

Our numerical computation depends, in a crucial way, on  the evaluation of a certain pairing in coherent cohomology on the mod $p$ fiber of a modular curve;
this evaluation is postponed to \S \ref{sec:Merel}, where we do it by relating it to Merel's remarkable computation \cite{Merel}. 
It is worth emphasizing how important Merel's computation is for us: it seemed almost impossible
to carry through our computation until we learned about Merel's results. Indeed, the role that Merel's computation 
plays here suggests that it would be worthwhile to understand how it might generalize to Hilbert modular surfaces. 
  
The expression of the derived Hecke algebra action as a cohomological cup product (see \eqref{above2}) strongly suggests a surprising relation with special values of the {\it $p$-adic triple product $L$-function}.  We hope to explore this relation in forthcoming work with Henri Darmon.

\subsection{Acknowledgements}
We would like to thank Brian Conrad for helpful explanations about flat cohomology,
  Frank Calegari  for pointing out Merel's paper and for interesting discussions related to that paper, and Shekhar Khare for several helpful and inspiring discussions
  (in particular related to \S \ref{GDS}).  We would also like to thank Nick Rozenblyum for help with Remark \ref{millenials}, and Henri Darmon for reading and providing comments
  on the paper.  Finally, we thank Loic Merel and Emmanuel Lecouturier  for exchanges concerning their work. 
  
\section{Notation} \label{notn}

  \subsection{Modular curves}
Fix an integer $N$. We will generically use the letter $R$ for a   $\Z[1/N]$-algebra.  For $R \rightarrow R'$ a morphism of $\Z[1/N]$-algebras and $Y$ an $R$-scheme, we denote by $Y_{R'}$
the base extension of $Y$ to $R'$.

Let $X = X_1(N)$ be the compactification (\cite{DR}) of the modular curve   parametrizing elliptic curves with with an $N$-torsion point.
We may construct $X$ as a smooth proper relative curve over $\Spec \ \Z[\frac{1}{N}]$, and the cusps give rise
 to a relative divisor $D \subset X$.   In particular, we obtain an $R$-scheme $X_R$ for each $\Z[1/N]$-algebra $R$. 
      We denote the universal generalized elliptic curve over $X$
by $\cA \longrightarrow X$.

  We will denote by 
  $$X_{01}(qN), X_1(qN)$$
  the modular curves that correspond to  adding $X_0(q)$ and $X_1(q)$ level structure to $X$.

Let $\omega$ be the line bundle over $X$ whose sections are given by weight $1$ forms. 
  More precisely, when $X$ is the modular curve, let $\Omega_{\cA/X}$ denote the relative cotangent bundle of $\cA/X$, pulled back to $X$ via the identity section.   We define $\omega$ as the pullback of $\Omega_{\cA/X}$ via the zero section $X \ra \cA$. 
  Therefore the sections of $\omega$ correspond to weight $1$ forms, whereas sections of 
    $\omega(-D)$ corresponds to weight $1$ {\em cusp} forms. 
Moreover, there is an isomorphism of line bundles  (\cite{K}, (1.5), A (1.3.17)):
\begin{equation} \label{spin structure}  \omega \otimes \omega(-D) \simeq \Omega^1_{X \rightarrow \Z[\frac{1}{N}]},\end{equation}
which says that  ``the product of a weight $1$ form and a cuspidal weight $1$ form
is a cusp form of weight $2$.''

  Let $\pi: X_R \rightarrow \Spec R$ be the structure morphism and consider the space of weight $1$ forms
  over $R$, in cohomological degree $i$ -- formally:
    $$ \Gamma(\Spec(R),\  R^{i} \pi_* \omega).$$
    We will denote this space, for short, by $H^i(X_R, \omega)$, and use similar notation for $\omega(-D)$. 
Therefore $H^0(X_R, \omega)$ (respectively $H^0(X_R, \omega(-D))$)  is the usual space   of weight $1$ modular forms (respectively cusp forms) with coefficients in $R$.

  \subsection{The residue pairing} \label{residuepairing}
 The pairing \eqref{spin structure} induces 
  $$ \pi_* \omega \otimes \mathrm{R}^1 \pi_* \omega(-D) \longrightarrow \mathrm{R}^1 \pi_* \Omega^1_{X_R/R}$$
  Since $X_R$ is a projective smooth curve over $R$, there is a canonical identification of the last factor with  the trivial line bundle; thus we get a pairing       
   $$H^0(X_R, \omega) \times H^1(X_R, \omega(-D)) \rightarrow  R,$$
  which we denote as $[- , -]_{\res,R}$.    (Here $\res$ stands for ``residue.'') 
  This pairing is compatible with change of ring, and if  $R$ is a field it is a perfect pairing.

  \subsection{The fixed weight one form $g$}    \label{gdef}

We want to localize our story throughout at a single weight one form $g$. 
  Therefore, fix
  $ g = \sum a_n q^n$
  a Hecke newform of level $N$ and Nebentypus $\chi$, normalized so that $a_1 = 1$. 
  Here $\chi$ is a Dirichlet character of level $N$. 
  
 We regard the $a_n$ as lying in some number field $E$, and indeed in the integer ring $\mathcal{O}$ of $E$.
 Thus $g$ extends to a section: 
 $$ g \in  H^0(X_{\mathcal{O}[\frac{1}{N}]}, \omega(-D)).$$ 
 
 We shall denote by 
 $H^*(X_{\mathcal{O}[\frac{1}{N}],} \omega)[g]$,
  that part of the cohomology which transforms under the Hecke operators in the same way as $g$,
  i.e. the common kernel of all $(T_{\ell}-a_{\ell})$ over all  primes $\ell$ not dividing $N$. 
  
 Extending $E$ if necessary,
 we may suppose that   
one can attach to $g$  a Galois representation, unramified away from $N$  \cite{DS}:
\begin{equation} \label{rhodef}  \rho : \Gal(L/\Q) \longrightarrow \GL_{2}(\mathcal{O})\end{equation}
  where $L$ is a Galois extension of $\Q$.  Here      the Frobenius trace  of $\rho$ at ${\ell}$ coincides with $a_{\ell}$, and 
the image of  complex conjugation $c$ under $\rho$ is conjugate to $\left( \begin{array}{cc} 1 & 0 \\ 0 & -1 \end{array}\right)$. 
(In the body of the text, we will primarily use the field cut out by the {\em adjoint} of $\rho$,
and  we could replace $L$ by this smaller field.)

  We emphasize the distinction between $E$ and $L$:
  \begin{quote}
  $E$ is a coefficient field for the weight one form $g$, and $L$ is the splitting
  field for the Galois representation of $g$.
  \end{quote}
  In our numerical examples, we will have $E=\Q$ and $L$ the Galois closure of a cubic field. 
  
    It will be convenient to denote 
 by $\Ad^0 \rho$ the trace-free adjoint of $\rho$, i.e. the associated action of $\Gal(L/\Q)$ on
 $2 \times 2$ matrices of trace zero and entries in $\mathcal{O}$. We denote by $\Ad^* \rho$ the $\mathcal{O}$-linear dual to $\Ad^0 \rho$, i.e.
 $$ \Ad^* \rho = \Hom(\Ad^0 \rho, \mathcal{O}),$$
a locally free $\mathcal{O}$-module endowed with an action of $\Gal(L/\Q)$.   
\footnote{We apologize for the perhaps pedantic distinction between $\Ad^* \rho$ and $\Ad^0 \rho$. Since we will
shortly be localizing at a prime larger than $2$, one could identify them by means of the
pairing $\mathrm{trace}(AB)$.  However, when working in a general setting,
one really needs to use $\Ad^*$, and following this convention makes it easier to compare with \cite{V3}.}
 
 For later use, it is convenient to choose a dual form $g'$ which will be paired with $g$ eventually.   In order
 that a Hecke equivariant pairing between $g$ and $g'$ be nonzero we should take
 $g'$  to be the form corresponding to the contragredient automorphic representation, i.e. 
  $$g' := \sum  \overline{a_n}  q^n  \in H^0(X_{\mathcal{O}[\frac{1}{N}]}, \omega(-D))$$
  where $\alpha \mapsto \overline{\alpha}$ is the complex conjugation in the CM field $E$. (In our examples,
  $E = \Q$, and therefore $g'=g$).

 \subsection{The prime $\mathfrak{p}$} \label{lambdaconditions}
 
   Let $\mathfrak{p}$ be a prime of $E$, 
above the rational prime $p$.  We make the following assumptions:
\begin{itemize}
\item  All weight one forms in characteristic $p$ lift to characteristic zero,
i.e. the natural map $$H^0(X_{\mathbb{Z}_p}, \omega) \rightarrow H^0(X_{\mathbb{F}_p},\omega)$$
is surjective. 
 
\item $p \geqslant 5$. 
\item $p$ is unramified inside $E$. 
\item There are no $p$th-roots of unity inside $L$. 
\item $p$ does not divide the order $[L:\Q]$. 
\end{itemize}

The representation $\rho$  may be reduced modulo $\mathfrak{p}$, obtaining  $\rhobar: G_{\Q} \longrightarrow \GL_2(\mathbf{F}_{\mathfrak{p}})$, 
 where $\mathbf{F}_{\mathfrak{p}} = \mathcal{O}/\mathfrak{p}$ is the residue field at $\mathfrak{p}$.    As before
 we may define the trace-free adjoint $\Ad^0 \rhobar$ and its dual $\Ad^* \rhobar$.

 \subsection{Taylor Wiles primes} \label{TWprimesec}
  A {\em Taylor--Wiles prime 
  $q$ of level $n \geq 1$} for $g$, or more precisely
  relative to the pair $(g, \mathfrak{p})$, will be, by definition:
  \begin{itemize}
  \item 
a   rational prime $q \equiv 1$ modulo $p^n$,  relatively prime to $N$; 
\item the data of $(\alpha, \beta)  \in \mathbf{F}_{\mathfrak{p}}$ with $\alpha \neq \beta$, 
 such that $\rhobar(\mathrm{Frob}_q)$ is conjugate to $\left( \begin{array}{cc} \alpha & 0 \\ 0 & \beta \end{array} \right)$.
 \end{itemize}
 
  Thus whenever we refer to Taylor--Wiles primes, we always regard the ordered pair $(\alpha, \beta)$ as part of the data: this amounts
 to an ordering of the eigenvalues of Frobenius.

If $\mathfrak{p}, n, q$ have been fixed,  where $q$ is a Taylor--Wiles prime of level $n$ and $\mathfrak{p}$ is a prime of $\mathcal{O}$ as above,  it is convenient to use the following shorthand notation:

\begin{itemize}
\item  Write  $k =  \mathcal{O}/ \mathfrak{p}^n$.    \item  
  Write 
 $(\Z/q)^*_p$ for the quotient
   of $(\Z/q)^*$ of size $p^n$, so that there is a noncanonical isomorphism $(\Z/q)_p^* \cong \Z/p^n$.
    \item 
Write $$ k\langle 1 \rangle = k \otimes (\Z/q)^*_{p}, \ \ k \langle -1 \rangle = \Hom((\Z/q)_{p}^*, k).$$
These are isomorphic as abelian groups to $k$, but not canonically so. 
\item Similarly for a  $\Z$-module $M$ we shall write
$$ M \langle n \rangle = M \otimes_{\Z} k \langle n \rangle.$$ 
 Thus, for example, $\mathbb{F}_p\langle 1 \rangle$ is canonically identified
with the quotient of $(\Z/q)^*$ of size $p$. 
 \end{itemize} 

These notations clearly depends on $\mathfrak{p}, n, q$; however  we do not explicitly indicate this dependence.  %

\subsection{The Stark unit group}

Let $\UnitL$ be the group of units of the integer ring of $L$. 

The key group of ``Stark units'' that we shall consider is  the following $\mathcal{O}$-module: 
\begin{align}  \label{dododo0} \OLrho  :=  \left( \UnitL \otimes_{\Z}  \Ad^* \rho   \right)^{G_{L/\Q}} = 
 \left( \UnitL \otimes_{\Z} \Hom_{\mathcal{O}}(\Ad^0 \rho, \mathcal{O}) \right)^{G_{L/\Q}} \\ 
 \label{dododo}  \stackrel{\sim}{\rightarrow} \Hom_{\mathcal{O}[G_{L/\Q}]}(\Ad^0 \rho, \UnitL \otimes \mathcal{O})\end{align}
 where $\Ad^0$ is the conjugation action of the Galois group on trace-free matrices in $M_2(\mathcal{O})$. 
 
 For instance,
 in the examples of modular forms attached to cubic fields,
the group $\UnitL$ will amount to (essentially) the unit group of that cubic field.

\begin{Lemma} \label{motrankone} 
$\OLrho \otimes_{\Z} \Q$ is an $E$-vector space of rank $1$ and   $\OLrho \otimes_{\Z} \Z_p$
is a free $\mathcal{O} \otimes_{\Z} \Z_p$-module of rank one. 
\end{Lemma}

\proof
Fix any embedding $\iota: E \hookrightarrow \C$; the $E$-dimension of $\OLrho \otimes_{\Z} \Q$
then coincides with the complex dimension of Galois invariants on $\UnitL \otimes_{\Q} (\Ad^* \rho)^{\iota}$.

In general, for any number field $L$, Galois over $\Q$,  and  any $\eta: G_{L/\Q} \rightarrow \GL_m(\C)$ without a trivial subrepresentation,  
the  dimension of $(\UnitL \otimes \eta)^{G_{L/\Q}}$ is known to be the dimension
of invariants for complex conjugation in $\eta$.       (In fact this is a straightforward consequence of the unit theorem.) It follows, therefore,
that the $E$-dimension of $\OLrho \otimes_{\Z} \Q$ is $1$ as claimed.

 The final claim now follows since our assumption on $p$ means that $\UnitL$ is free of $p$-torsion. 
\qed

Fix a nonzero element $u \in \OLrho$ 
  in such a way that $[\OLrho: \mathcal{O} . u]$ is relatively prime to $p$ (e.g. a generator,
  if $\OLrho$ is a free $\mathcal{O}$-module). 
   Later we will also work with the  $\mathcal{O}$-dual 
\begin{equation} \label{Edual} 
\OLrho^{\vee} := \Hom_{\mathcal{O}}(\OLrho, \mathcal{O}), 
 \end{equation}
and denote by $u^* \in \OLrho^{\vee}$ a nonzero element, chosen so that $\langle u, u^* \rangle \in \mathcal{O}$
is not divisible by any prime above $p$.

\subsection{Aside: comparison of $\OLrho$ with the motivic cohomology group from \cite{V3}}   
This section is not used in the remainder of this paper. 
It serves to connect the previous construction with the discussion in   \cite{V3}:

  We may construct a $3$ dimensional Chow motive $\Ad^0 \Mg$, with coefficients in $E$, attached to the trace-free adjoint $\Ad^* \rho_g$ -- in other words,
  the {\'e}tale cohomology of $\Ad^0 \Mg$ is concentrated in degree zero and identified, as a Galois representation, with $\Ad^* \rho_g$. 
  
Now consider the motivic cohomology $H^1_{\mot}(\Q, \Mg(1))$,
or more precisely the subspace of integral classes $(-)_{\Int}$ described by Scholl \cite{Scholl}.

The general conjectures of \cite{V3, PV}, transposed to the current (coherent) situation, predict that the dual of  $H^1_{\mot}(\Q, \Mg(1))_{\Int}$  should act 
on  $H^*(X_E, \omega)[g]$.

There is a natural map 
\begin{align} H^1_{\mot}(\Q, \Mg(1))_{\Int} & \longrightarrow      H^1_{\mot}(L, \Mg(1))_{\Int}^{G_{L/\Q}}     \\ &=   \left( \UnitL \otimes \Ad^* \rho  \otimes \Q\right)^{G_{L/\Q}} = \OLrho \otimes \Q.
\end{align}

Although we did not check it, this map is presumably an isomorphism.  In the present paper we will never
directly refer to the motivic cohomology group. Rather 
we work with the right-hand side
(or its integral form $\OLrho$) as a concrete substitute for the motivic cohomology group.

  \subsection{Reduction of a Stark unit at a Taylor--Wiles prime $q$}

Let $q$ be a Taylor--Wiles prime (as in \S \ref{TWprimesec});  we shall define a canonical reduction map  $$ \theta_q:  \OLrho \longrightarrow  k \langle 1 \rangle.$$ 
  For example,
in the examples of modular forms attached to cubic fields,
this will amount to the reduction of a unit in the cubic field at a degree one  prime above $q$. 
Although explicit, the general definition is  unfortunately opaque (the motivation comes from computations in \cite{V3}).

For any prime $\mathfrak{q}$ of $L$ above $q$, with associated Frobenius element $\Frob_{\mathfrak{q}}$, let
  $D_{\mathfrak{q}}= \langle \Frob_{\mathfrak{q}} \rangle \subset \Gal(L/\Q)$ be the associated decomposition group, the stabilizer of $\mathfrak{q}$. 
We may construct  a  $D_{\mathfrak{q}}$-invariant element
\begin{equation} \label{eqdef} \mathsf{e}_{\mathfrak{q}}  = 2 \rho(\Frob_{\mathfrak{q}}) -  \mathrm{trace}\ \rho(\Frob_{\mathfrak{q}})  \in \Ad^0 \rho\end{equation}
where we regard the middle quantity as a $2 \times 2$ matrix with coefficients in $\mathcal{O}$ and trace zero, thus belonging to
$  \Ad^0  \rho.$ Pairing with $\mathsf{e}_{\mathfrak{q}}$ and reduction mod $\mathfrak{p}^n$ induces 
$$\mathsf{e}_{\mathfrak{q}}:  \Ad^*  \rho \longrightarrow \mathcal{O} \rightarrow k,$$
equivariantly for the Galois group of $\Q_q$.  Also, for $g \in \Gal(L/\Q)$ we have 
\begin{equation} \label{compati}   \mathsf{e}_{g \mathfrak{q}} = \Ad(\rho(g)) e_{\mathfrak{q}}.\end{equation}

Write $L_q = (L \otimes \Q_q)$ and let $\mathcal{O}_{L_q}$ be the integer subring thereof. 
Thus $\mathcal{O}_{L_q}/q \simeq \prod_{\mathfrak{q}|q} \mathbf{F}_{\mathfrak{q}}$. 
Fix a prime $\mathfrak{q}_0$ of $L$ above $q$. 
The inclusion of units for  the number field $L$ into local units $\mathcal{O}^*_{L_q}$ induces 
\begin{equation}\label{twmap0} (\OLrho) \rightarrow (  \prod_{\mathfrak{q}|q} \mathbf{F}_{\mathfrak{q}}^*  \otimes  \Ad^*\rho)^{G_{L/\Q}}
\stackrel{\sim}{\rightarrow}  \left(\mathbf{F}_{\mathfrak{q}_0}^* \otimes \Ad^* \rho \right)^{D_{\mathfrak{q}_0}} \stackrel{\mathsf{e}_{\mathfrak{q}_0}}{\longrightarrow}  (\mathbf{F}_{\mathfrak{q}_0}^* \otimes  k)^{D_{\mathfrak{q}_0}}  \rightarrow   k\langle1 \rangle. \end{equation}
where  the second map is projection onto  the factor corresponding to $\mathfrak{q}_0$.      

 The resulting composite is independent of the choice of $\mathfrak{q}_0$, because of \eqref{compati}. 
We call it $\theta$, or $\theta_q$ when we want to emphasize the dependence on the Taylor-Wiles prime $q$:
$$ \theta \mbox{ or } \theta_q: \OLrho \rightarrow k \langle 1 \rangle .$$

 \section{Derived Hecke operators and the main conjecture} \label{DHmain}

  We follow the notation of \S \ref{notn}; in particular,

\begin{itemize}
\item[-] $g$ is a modular form with coefficients in the integer ring $\mathcal{O}$;
we have associated to it a $\mathcal{O}$-module $\OLrho$ of ``Stark units'' of rank $1$.
 \item[-] Fixing a prime $\mathfrak{p}$ of $\mathcal{O}$, we will work with the coefficient ring   $k = \mathcal{O}/\mathfrak{p}^n$
 with residue field $\mathbb{F}_{\mathfrak{p}}$ of characteristic $p$. 
 \end{itemize}

 In this section we define derived Hecke operators and formulate the main conjecture 
 concerning their relationship to $\OLrho$.   This discussion is obtained by transcribing the theory of \cite{V3} to the present context;
 in this section, we just describe the conclusions of this process.

 For each $q \equiv 1$ modulo $p^n$
  and each $z \in k \langle -1 \rangle$, we will 
  produce an operator
$$ T_{q,z} : H^0(X_k, \omega) \rightarrow H^1(X_k, \omega). $$ 

Note that $q$ need not be a Taylor--Wiles prime  (in the sense of \S \ref{lambdaconditions}) for the definition of $T_{q,z}$ -- in other words, we do not use the assumption on the Frobenius element.
However, our conjecture pins down the action of $T_{q,z}$ only at Taylor--Wiles primes. 

\subsection{The Shimura class}  \label{Shimuraclassdef}

  Start with the Shimura covering $X_1(q) \rightarrow X_0(q)$,
  and pass to the unique subcovering with Galois group $(\Z/q)^*_p$;
  call this $ X_1(q)^{\Delta} \rightarrow X_0(q)$. 
    By  Corollary 2.3 of \cite[Chapter 2]{mazur},
  it extends to an {\'e}tale covering of schemes over $\Z[\frac{1}{qN}]$, and in particular    induces an {\'e}tale cover $ X_1(q)^{\Delta}_{k} \rightarrow X_0(q)_k$. 
  It therefore gives rise to 
 a class in the etale $H^1$, i.e.
   $$\Shim \in H^1_{\et}(X_0(q)_k, k \langle 1 \rangle)$$ 
In the category of {\'e}tale sheaves over  $X_0(q)_k$ there is a natural map $k \rightarrow \mathbb{G}_a$.   %
Then a class in $H^1_{\et}(X_0(q)_k, k \langle 1 \rangle)$ defines a class in 
$$H^1_{\et}(X_0(q)_k, \mathbb{G}_a \langle 1 \rangle) \simeq H^1_{\mathrm{Zar}}(X_0(q)_k, \mathcal{O}\langle 1 \rangle),$$
because of the coincidence of the {\'e}tale and Zariski cohomologies with coefficients in a quasi-coherent sheaf. 
This construction has thus given  a class, associated to the Shimura cover, but now in Zariski cohomology:
$$\Shim \in H^1_{\mathrm{Zar}}(X_0(q)_k, \mathcal{O}\langle 1 \rangle),$$
which we shall sometimes call the Shimura class. 

It is reassuring to note that $\Shim$ is in fact nonzero, even modulo the maximal ideal $\mathfrak{p}$,
as one sees by computing with the Artin--Schreier sequence
$$0 \rightarrow \mathbb{F}_{\mathfrak{p}} \rightarrow \mathbb{G}_a \rightarrow \mathbb{G}_a \rightarrow 0$$
over $k/(\mathfrak{p})$.

 \subsection{Construction of the derived Hecke operator}

 The class $\Shim$ just defined  can be pulled back to $H^1_{\mathrm{Zar}}(X_{01}(qN)_k, \mathcal{O} \langle 1 \rangle)$. 
 We denote this class by $\Shim_X$, to distinguish it from $\Shim$ at level $q$. 
 
  Thus cup product with $\Shim_X$ gives a mapping
 $$H^0(X_{01}(qN)_k, \omega) \stackrel{\cup \Shim_X}{\longrightarrow} H^1(X_{01}(qN)_k, \omega \langle 1 \rangle).$$
 Finally, to obtain the derived Hecke operator we 
add a push-pull as in the usual Hecke operator definition: 
\begin{equation} \label{dhdefn} H^0(X_k,\omega) \stackrel{\pi_1^*}{\longrightarrow} H^0(X_{01}(qN)_k, \omega)\stackrel{\cup \Shim_X}{\longrightarrow} H^1(X_{01}(qN)_k, \omega) \langle 1 \rangle 
\stackrel{\pi_{2*}}{ \longrightarrow} H^1(X_k, \omega) \langle 1 \rangle,\end{equation}
where $\pi_1, \pi_2: X_{01}(qN) \rightarrow X$ are the two natural degeneracy maps (at the level
of the upper half-plane, we understand $\pi_1$ to be $z \mapsto z$, and $\pi_2$ to be $z \mapsto qz$.)
  Observe that without the middle $\cup \Shim_X$ this would be the usual Hecke operator at $q$.  
In other words, we have constructed a map \begin{equation} \label{first} H^0(X_k, \omega)  \rightarrow H^1(X_k, \omega) \langle 1 \rangle,\end{equation}
 and correspondingly for  $z \in k \langle -1 \rangle$ we will denote by $T_{q,z}$ the corresponding ``derived Hecke operator''
 \begin{equation} \label{dh2} T_{q,z} : H^0(X_k, \omega) \rightarrow  H^1(X_k, \omega). \end{equation} 
 obtained by multiplying \eqref{first} by $z$. 
 
  Although by presenting the bare definition the construction may seem a little {\em ad hoc}, this definition
  is really a specialization of the general theory of \cite{V3}, and is indeed very natural. We explain this in more detail in \S \ref{final}.

  \subsection{The conjecture}
 We now formulate the main conjecture.  It asserts that the various operators $T_{q,z}$ all fit together into a single
 action of $\OLrho^{\vee}$ on the $g$-part of cohomology. As formulated in \cite{V3},
 the conjecture is ambiguous up to a rational factor, and we will not attempt
 to remove this ambiguity here (although our computations suggest that this factor might  have a simple description). 

Terminology:
\begin{itemize}
\item Suppose that $\alpha \in E$
 and $V$ is a  $k$-module. For $x,y \in V$, we will write
\begin{equation}  \label{prop test} x = \alpha y \end{equation}
 if we may write $\alpha = A/B$, where $A, B \in \mathcal{O}$ 
 are not both divisible by $\mathfrak{p}$,
 in such a way that $\bar{B} x=\bar{A} y$ (here $\bar{A}, \bar{B}$ are the reductions of $A, B$
 under $\mathcal{O} \rightarrow k$.) 
 
 In particular, if $V$ is a $k$-line this has the following meaning:
 \begin{itemize}
 \item if $x=y=0$, then \eqref{prop test} is understood to always be true. 
 \item Otherwise, we can make sense of $[x:y] \in \mathbb{P}^1(k)$, and \eqref{prop test}
 means that the reduction of $\alpha \in \mathbb{P}^1(E) \rightarrow \mathbb{P}^1(k)$
 equals $[x:y]$. 
 \end{itemize}
 
 \item For $h \in H^*(X_{\mathcal{O}[\frac{1}{N}]}, \omega)$
 we write $\bar{h}$ for the reduction of $h$ to $H^*(X_k, \omega)$. 
 
 \item Recall that we defined a  reduction map $\theta_q: \OLrho \rightarrow k \langle1 \rangle$.  Also, 
the pairing between $\OLrho$ and $\OLrho^{\vee}$, which is perfect after localization at $p$, descends to a perfect pairing
on $\OLrho \otimes k$ and $\OLrho^{\vee} \otimes k$. 
With respect to this pairing, the map $\theta_q$ has an adjoint:
$$\theta_q^{\vee}:  k \langle -1 \rangle \rightarrow \OLrho^{\vee} \otimes k.$$
 Explicitly  for $z \in k\langle -1 \rangle$, 
\begin{equation} \label{theta vee explicit} \theta_q^{\vee}(z) =  u^* \otimes  \frac{  \langle z, \theta_q(u) \rangle}{\langle u^*, u \rangle},\end{equation}
where $u \in \OLrho,  u^* \in \OLrho^{\vee}$ are as defined around \eqref{Edual}. 
 \end{itemize}

   \begin{Conjecture}\label{Mainconj} There is an action $\star$ of $\OLrho^{\vee}$
   on $H^*(X_{\mathcal{O}[\frac{1}{N}],} \omega)[g]$,  and $\alpha \in E$ such that for every     $(\mathfrak{p}, n, q, z)$, with
   \begin{itemize}   \item[-] $\mathfrak{p}$   a prime of $E$ satisfying the conditions of \S \ref{lambdaconditions};
   \item[-] $n \geqslant 1$  an integer; 
   \item[-]  $q$  a Taylor--Wiles prime of level $n$, in particular $q \equiv 1 (p^n)$. 
\item[-]   $z \in (\mathcal{O}/\mathfrak{p}^n) \langle -1 \rangle$,
\end{itemize}
we have the following equality:
\begin{equation} \label{equa}   T_{q,z} \bar{g}   = \alpha  \overline{\left( \theta^{\vee}_q(z)^{\sim} \star g \right)},\end{equation} 
 
 On the right hand side, $\theta^{\vee}_q(z)^{\sim}$ means that we choose an arbitrary lift of $\theta^{\vee}_q(z) \in \OLrho^{\vee} \otimes k$
to $\OLrho^{\vee}$, and the bar refers to reduction mod $\mathfrak{p}^n$.

   \end{Conjecture}

 In what follows, we will write \eqref{equa} in the abridged form
\begin{equation}  \label{proptoequality}T_{q,z} g \propto   \overline{ \theta^{\vee}_q(z)  \star g}.\end{equation}
  The meaning here is 
  that equality holds, in the sense described above,
for some fixed coefficient of proportionality $\alpha \in E$.  (Note we have suppressed explicit mention of the lift $\theta^{\vee}_q(z)^{\sim}$ from the notation;
in any case the right-hand side is independent of this choice of lift.)

\section{Relationship to Galois deformation theory} \label{GDS}
 In this section -- {\em which is not used in the rest of the paper} -- we  shall sketch  a proof that,  
 in the case $n=1$, 
\begin{equation} \label{rq} \mbox{ vanishing of $T_{q,z} \bar{g}$} \implies
 \mbox{ vanishing of $\theta_q: \OLrho \rightarrow k\langle 1 \rangle$. } 
 \end{equation} 
 assuming an ``$R=T$'' theorem for weight one forms at the level of $g$, as well as further technical conditions. Such a theorem
 is known in some generality by the work of Calegari \cite{Calegari}.

     This result (and its proof) is in line with results and proofs  from \cite{V3}. Indeed,  our methods would show that \eqref{rq} is an equivalence,  {\em if} we knew an  ``$R=T$'' theorem for weight one forms 
with (Taylor--Wiles) auxiliary level. 

\subsection{Setup}

Let $q$ be a prime such that the eigenvalues of $\rhobar$ on the Frobenius at $q$ are distinct elements of $\mathbf{F}_{\mathfrak{p}}$, 
 say $\alpha$ and $\beta$. 
   Let $\mathfrak{m}$ be the ideal of the Hecke algebra associated to the Galois representation $\rhobar$. 

In addition to the conditions from \eqref{lambdaconditions}
we assume that:

\begin{itemize}
\item[(i)] $n=1$ so that $k = \mathcal{O}/\mathfrak{p}$ is a field. 
\item[(ii)]  For each prime $\nu$ dividing $N$,   the residual representation $\rhobar$ is of the form
$\chi_1 \oplus\chi_2$, where $\chi_1$ is ramified and $\chi_2$ is unramified. 
\item[(iii)] $p$ does not divide $\nu-1$, for each $\nu$ as above. 
\item[(iv)] $p$ does not divide $[L:\Q]$, and does not divide the order of the class group of $L$. 
\item[(v)]   The $\mathfrak{m}$-completion of the space of modular forms at level $\Gamma_1(N)$,  with coefficients in $\mathcal{O}$, 
is free rank one over $\mathcal{O}_{\mathfrak{p}}$.   (In particular, there are no congruences modulo $\mathfrak{p}$ between $g$ and other
weight one forms, either in characteristic zero or characteristic $p$.) 
\end{itemize}

   Let $\mathfrak{m}_{\alpha}$ be the maximal ideal of the Hecke algebra
for $X_{01}(qN)$  obtained by adjoining $U_q - \alpha$ to the ideal $\mathfrak{m}$;
similarly we define $\mathfrak{m}_{\beta}$. These ideals also have evident analogues
where we add $\Gamma_1(q)$ level to $X$, rather than just $\Gamma_0(q)$ level, and we denote these analogues by the same letters.  
   
Our assumption (v), and the assumption of torsion-freeness from \eqref{lambdaconditions},  means that  {\small \begin{equation} \label{assumpt} \dim H^0(X_k, \omega)_{\mathfrak{m}} = \dim H^0(X_{01}(qN)_k, \omega)_{\mathfrak{m}_{\alpha}} =   
\dim H^0(X_{01}(qN)_k, \omega)_{\mathfrak{m}_{\beta}}  = 1,\end{equation}}
i.e. all three spaces above are $k$-lines; the same statement is true for $H^1(-)$. 

Let $g_{\alpha}$ and $g_{\beta}$, respectively, span the second and third spaces in the line above. 
Therefore $U_q g_{\alpha} = \alpha g_{\alpha}$ and $U_q g_{\beta} = \beta g_{\beta}$. 
We normalize these so that $\pi_1^* g = g_{\alpha} + g_{\beta}$.

 Since the pushforward $\pi_{1*}$  via the natural projection $\pi_1: X_{01}(qN) \rightarrow X$ induces an isomorphism on  each of the $U_q$-eigenspaces, $g_{\alpha} \cup \Shim_X$ vanishes
if and only if $\pi_{1*}(g_{\alpha} \cup \Shim_X)$ vanishes.  Observe $(U_q - \beta) \pi_1^* g = (\alpha -\beta) g_{\alpha}$. We are assuming $\alpha \neq \beta$  and therefore
\begin{equation} \label{bre}  g_{\alpha} \cup \Shim_X = 0 \iff \pi_{1*} ( (U_q - \beta) \pi_{1}^* g \cup \Shim_X)  =0.\end{equation}
 Now $\pi_{1*} (\pi_1^* g \cup \Shim_X) = g \cup (\pi_{1*} \Shim_X) = 0$ and $\pi_{1*} \Shim_X$ is trivial:

 \begin{lemma} The pushforward of $\Shim_X$
 by the natural projection $\pi: X_{01}(qN) \rightarrow X$ is trivial.
\end{lemma}
  \proof  The existence of the trace  map  \cite[Expose 17, Section 6.2]{SGA4}
   gives a map $\pi_* (\Z/p) \rightarrow \Z/p$ of {\'e}tale sheaves, 
 compatible with the usual trace $\pi_* \mathbb{G}_a \rightarrow \mathbb{G}_a$. 
 For this reason, it is sufficient to show that the (trace-induced) map
 $$H^1_{\et}(X_{01}(qN)_k, (\Z/q)^*_p) \rightarrow H^1_{\et}(X_k, (\Z/q)^*_p)$$
 pushes the Shimura class forward to the trivial class.

 If $\iota$ is the inclusion of  an open curve into a complete curve induces  then $\iota^*$ is an injection on $H^1$.  Therefore, it suffices to
 show a similar statement for the open modular curves; restricted to these, the map $\pi$ is {\'e}tale.

  Define finite groups 
 $$G = \GL_2(\Z/q\Z) \supset  B = \left(\begin{array}{cc} * & * \\ 0 & * \end{array} \right), $$
Then $X_{01}(qN)$ and $X_k$ are quotients of a suitable modular curve
by $B$ and $G$ respectively.
 This allows to reduce to verifying the triviality of the transfer  in group cohomology, from $B$ to $G$, of 
$\alpha   \in H^1(B, (\Z/q)_p^*)$, defined via  
$$ \alpha: \left(\begin{array}{cc} a & c \\ 0 & d \end{array} \right) \mapsto (a/d) \in (\Z/q)^*_p.$$
 This is a straightforward computation. \qed

Continuing from \eqref{bre}, we find
\begin{equation} \label{bre2} g_{\alpha} \cup \Shim_X = 0 \iff \pi_{1*} (U_q \pi_1^* g \cup \Shim_X) = 0.\end{equation}
The final expression can be verified to be an invertible multiple of $T_{q,z} g$ for 
some nonvanishing $z \in \mathcal{O}/\mathfrak{p} \langle -1 \rangle$.   Therefore:
\begin{equation} \label{bre3} g_{\alpha} \cup \Shim_X = 0 \iff T_{q,z} \bar{g} = 0.\end{equation}

  Write $\Delta =  (\Z/q)^*_p$; since we are assuming that $n=1$,
  the group $\Delta$ is cyclic of order $p$ and we have an isomorphism
  $k[\Delta] \simeq k[T]/T^p$, whose inverse sends $T$ to $\delta -1$, for any generator $\delta$ of $\Delta$. 
  Let $X_1(Nq)^{\Delta}$ be the subcovering of $X_1(Nq) \rightarrow X_{01}(qN)$
  which corresponds to the quotient $(\Z/q)^* \rightarrow (\Z/q)^*_p$
  of deck transformation groups.

  \begin{lemma} \label{cupproductlemma}
  The cup product $\cup \Shim_X$
is nonzero as a map on $H^*(X_{01}(qN)_k, \omega)_{\mathfrak{m}_{\alpha}}$ if and only if
\begin{equation} \label{diamond} \dim H^0(X_1(Nq)_k^{\Delta}, \omega)_{\mathfrak{m}_{\alpha}} =1.\end{equation}
  \end{lemma}

The usual Taylor-Wiles method, for classical modular forms on $\GL_2$,
relies crucially on producing ``more'' modular forms when adding ``$\Gamma_1(q)^{\Delta}$ level'' at auxiliary primes $q$.
Thus the Lemma says: the derived Hecke operator is nontrivial precisely when this {\em fails}, a failure  that is rectified
in the Calegari--Geraghty approach \cite{CG}. 
  \proof    
By  the methods of \cite{CG}, we may obtain   
   a complex $C$ of  free $k[\Delta]$-modules (with degree-decreasing differential)
    together with isomorphisms:
\begin{equation} \label{r4} H^i \Hom_{k[\Delta]}(C, k) \simeq   H^i(X_{01}(qN)_k, \omega)_{\mathfrak{m}_{\alpha}},\end{equation}
 With  reference to the latter isomorphism, 
cup product with $\Shim_X$ on the right is represented
by the natural action of a nontrivial class in $\Ext^1_{k[\Delta]}(k, k)$ on the left hand side
(note that $H^i \Hom_{k[\Delta]}(C,k)$ is identified with homomorphisms
from $C$ to $k[i]$ in the derived category.)

 Replacing $C$ by a minimal free resolution  we may assume that $C$
 is the complex given by 
 $$ k[\Delta] \stackrel{A}{\leftarrow} k[\Delta],$$ where $A \in k[\Delta]$ belongs to the augmentation ideal. 
  Under the identification of $k[\Delta]$ with $k[T]/T^p$,
 the element $A$ corresponds to an invertible multiple of $T^i$, for some $0 \leq i \leq p-1$, and then \eqref{r4} implies
  $$ \dim H^0(X_1(Nq)_k^{\Delta}, \omega)_{\mathfrak{m}_{\alpha}} = i. $$ 
    
We shall show that cup product with $\Shim_X$ is nontrivial if and only if $i = 1$. 
To compute the action of $\Ext^1_{k[\Delta]}(k,k)$ we may consider the following diagram:  

    \begin{equation}
  \xymatrix{
 C  &  &  0   \ar[d]  &  k[\Delta]     \ar[l]   \ar[d]^1   & k[\Delta] \ar[d]^{T^{i-1}} \ar[l]^{T^i} &  \ar[l]  \ar[d] 0      \\
  k   &  &  0  \ar[d]      & k[\Delta]  \ar[l]  \ar[d] & k[\Delta] \ar[l]^{T}  \ar[d]^{1}  & \ar[l]^{T^{p-1}} k[\Delta]  \ar[d]^{T^{p-2}} & \ar[l]^{T} k[\Delta]  \ar[d]^1& \ar[l]  \dots \\
 k[1]  &  &   0   & \ar[l] 0  &  k[\Delta] \ar[l]  & k[\Delta] \ar[l]^{T} & \ar[l]^{T^{p-1}} k[\Delta]  & \ar[l] \dots 
  }
  \end{equation}
The horizontal complexes are, respectively, $C$, a projective resolution of $k$, and a projective resolution of $k[1]$. 
Continuing to take $\Hom$ in the derived category of $k[\Delta]$-modules, the top vertical map of complexes represents a generator for  $\Hom(C, k)$
  and the bottom vertical map of complexes represents a nontrivial class in $\Ext^1_{k[\Delta]}(k,k) = \Hom(k, k[1])$.
Therefore, the composite map in $\Hom(C, k[1])$ is represented by the diagram  
    \begin{equation}
  \xymatrix{
 &  0   \ar[d]  &  k[\Delta]     \ar[l]   \ar[d]  & k[\Delta] \ar[d]^{T^{i-1}} \ar[l]^{T^i} &  \ar[l]  \ar[d] 0      \\
  &   0   & \ar[l] 0  &  k[\Delta] \ar[l]  & k[\Delta] \ar[l]^{T}
  }
  \end{equation}
This is nullhomotopic exactly when $T^{i-1}$ is divisible by $T$, i.e. $i \geq 2$. \qed

Taking \eqref{bre3} together with the Lemma, we see
$$ T_{q,z} \bar{g} \neq 0 \iff \dim H^0(X_1(Nq)_k^{\Delta}, \omega)_{\mathfrak{m}_{\alpha}} =1.$$

Consider the map
$$ f: \mathrm{R}'  \otimes k \rightarrow \mathrm{R} \otimes k,$$
  where $\mathrm{R}$ (resp. $\mathrm{R}'$) are the weight 1, determinant $\chi$, deformation rings for $\rhobar$
at level $\Gamma_1(N)$ and with level $\Gamma_1(Nq)$ respectively. The local conditions $\mathcal{L}$ for $R$
and $\mathcal{L}'$ for $R'$ are as follows:
\begin{itemize}
\item  At  $p$, we require that deformation remains unramified. 
\item At $q$ we impose unramified  for $R$ and no condition for $R'$. 
\item 
For primes $\nu$ dividing $N$, we do not need to impose any condition: We have assumed that $\rhobar$ is a direct sum $\chi_1 \oplus \chi_2$ of two characters,
with $\chi_1$ ramified and $\chi_2$ unramified. In particular 
$ H^1(\Q_{\nu}, \Ad^0 \rhobar)$ is $1$-dimensional, corresponding to deforming
 $\chi_1 \leftarrow \chi_1 \psi, \chi_2 \leftarrow \chi_2 \psi^{-1}$ for a character $\psi$ with trivial reduction.  In \cite{Calegari}
 the assumption is imposed that in fact $\chi_2  \psi^{-1}$ remains unramified, but we do not 
 need to explicitly impose this because  we assumed that  $p$ is relatively prime to $\nu-1$ -- thus  the character $\psi$ is  automatically unramified at $\nu$. 
 In particular, we have automatically
 $$H^1(\Q_{\nu}, \Ad^0 \rhobar) = H^1_{\ur}(\Q_{\nu}, \Ad^0 \rhobar),$$
 where we recall that
 for a module $M$ under the Galois group of $\Q_{\ell}$,
 the ``unramified'' classes $H^1_{\ur} \subset H^1$ are defined to be those that arise from inflation from the Galois cohomology of $\F_{\ell}$ acting on inertial invariants on $M$. 
  
 \end{itemize}

 Assuming an $R=T$ theorem
for $g$, we  have $R \otimes k = k$.   
The map on tangent spaces induced by $f$, call it $f^*$,  fits into the following diagram, 
with reference to the usual identification of  tangent spaces with  Galois cohomology:  
{\small  \begin{equation*}
\begin{split} H^1_{\mathcal{L}}(\Q, \Ad^0 \rhobar) \stackrel{f^*}{\hookrightarrow } H^1_{\mathcal{L}'}(\Q, \Ad^0 \rhobar)
 \rightarrow \frac{H^1(\mathbb{Q}_q, \Ad^0 \rhobar)}{H^1(\Z_q, \Ad^0 \rhobar)} \stackrel{j}{\rightarrow} \\ H^2_{\mathcal{L}}(\Q, \Ad^0 \rhobar)
 \rightarrow H^2_{\mathcal{L}'}(\Q, \Ad^0 \rhobar)\end{split} \end{equation*}}

 $f^*$ is surjective exactly when $j$ is injective. Since the middle group in the exact sequence is one-dimensional, injectivity of $j$ is the same as
 nonvanishing of $j$.  Under Tate global duality, the map $j$ is dual to
\begin{equation} \label{dual reduction} H^1_{\mathcal{L}'^{\vee}}(\Q, \Ad^* \rhobar(1))  \stackrel{j^{\vee}}{ \rightarrow} H^1(\mathbf{F}_q, \Ad^* \rhobar(1)),\end{equation}
 where $\mathcal{L}'^{\vee}$ is the dual condition to $\mathcal{L}'$: it refers to classes that are unramified at primes not dividing $pN$,  unramified (equivalently: trivial) at primes dividing $N$,   and at $p$  belong to the Bloch-Kato $f$-cohomology (a more concrete description is given below).

We will show in the next subsection that:
\begin{equation} \label{altst}\mbox{ the map $j^{\vee}$ vanishes exactly when 
$\theta_q: \OLrho \rightarrow k\langle 1 \rangle$ does.}
\end{equation}
  Therefore, the nonvanishing of $\theta_q$ implies the
injectivity of $j$, which implies the 
 surjectivity of $f^*$,
which implies $R' \otimes k = k$, which implies  \eqref{diamond} by a multiplicity one argument. Then \eqref{bre3} and Lemma \ref{cupproductlemma} show that $T_{q,z} \bar{g} \neq 0$ as desired.

That concludes our  proof for \eqref{rq}; 
note finally that if we had a theorem $R'=T'$ all this reasoning would be reversible and we get
an equivalence in \eqref{rq}. 
\subsection{Relation of $\OLrho$ to Galois cohomology}
To conclude we must relate $j^{\vee}$ and $\theta_q$, and thereby prove \eqref{altst}.

 As in \eqref{rhodef}, $\rho$ is a representation into $\GL_2(\mathcal{O})$;
 let $\rho_{\mathfrak{p}}$ be the same representation, but now considered as valued in $\GL_2(\mathcal{O}_{\mathfrak{p}})$; thus
 $$ \Ad^* \rho_{\mathfrak{p}} = \Ad^* \rho \otimes_{\mathcal{O}} \mathcal{O}_{\mathfrak{p}}.$$
 We write $(\OLrho)_{\mathfrak{p}}$ for $\OLrho \otimes_{\mathcal{O}} \mathcal{O}_{\mathfrak{p}}$.

Consider
 $$ H^1_{\mathrm{ur}}(\Q, \Ad^* \rho_{\mathfrak{p}}(1))  $$ 
 where the subscript $\ur$ means that we consider classes that are unramified at primes away from $p$,
and, at $p$, belong the Bloch-Kato $f$-space.  (What this means is made explicit in the computation of $H^1_{\ur}(L, \mathcal{O}_{\mathfrak{p}}(1))$ in the diagram below.)

Restriction to $L$ gives horizontal maps in the following diagram
$$%
\xymatrix{
0 \ar[r] \ar[d]  & 0 \ar[d] \\ 
H^1_{\ur}(\Q, \Ad^* \rho_{\mathfrak{p}}(1)) \ar[d]  \ar[r]^{i  \qquad \qquad \qquad} &   \ar[d] [\Ad^* \rho_{\mathfrak{p}} \otimes_{\mathcal{O}_{\mathfrak{p}}}  H^1_{\ur}(L, \mathcal{O}_{\mathfrak{p}}(1))]^{G_{L/\Q}} =  (\OLrho)_{\mathfrak{p}}   \\ 
H^1(\Q, \Ad^* \rho_{\mathfrak{p}}(1))   \ar[d] \ar[r]^{j \qquad \qquad} &   \ar[d]  [\Ad^* \rho_{\mathfrak{p}} \otimes _{\mathcal{O}_{\mathfrak{p}}} H^1(L, \mathcal{O}_{\mathfrak{p}}(1))]^{G_{L/\Q}}   \\
\bigoplus_v \frac{H^1(\Q_v, \Ad^* \rho_{\mathfrak{p}}(1))}{H^1_{\ur}(\Q_v, \Ad^* \rho_{\mathfrak{p}}(1))}  \ar[r]^k &   \left[ \bigoplus_w \frac{H^1(L_w, \Ad^* \rho_{\mathfrak{p}}(1))}{H^1_{\ur}(L_w, \Ad^* \rho_{\mathfrak{p}}(1))}  \  \right]^{G_{L/\Q}} . 
}$$%
The vertical columns are exact at top and middle, and, in the bottom row, the sum is taken over all places $v$ of $\Q$,
and then over  all places $w$ of $L$. 
 
\begin{lemma}
$i$ induces an isomorphism $ H^1_{\mathrm{ur}}(\Q, \Ad^* \rho_{\mathfrak{p}}(1))  \simeq (\OLrho)_{\mathfrak{p}}$.
Also, as long as the class group of $L$ is prime to $p$, the reduction  modulo $\mathfrak{p}$ map 
 $ H^1_{\mathrm{ur}}(\Q, \Ad^* \rho_{\mathfrak{p}}(1))   \rightarrow H^1_{\mathrm{ur}}(\Q, \Ad^* \rhobar(1)) $
 is surjective. \end{lemma}

\proof 
The map $j$  is an isomorphism by considering the inflation-restriction sequence: the group $G_{L/\Q}$
has order prime to $p$.  

This means $i$ is injective.
 $i$ will be surjective if $k$ is injective. In fact, for a place $\mathfrak{q}$ of $L$ above $v$, the map
\begin{equation} \label{woo-hoo} \frac{H^1(\Q_v, \Ad^* \rho_{\mathfrak{p}}(1))}{H^1_{\ur}(\Q_v, \Ad^* \rho_{\mathfrak{p}}(1))} \rightarrow \left[ \frac{H^1(L_{\mathfrak{q}}, \Ad^* \rho_{\mathfrak{p}}(1))}{H^1_{\ur}(L_{\mathfrak{q}}, \Ad^* \rho_{\mathfrak{p}}(1))}  \right].\end{equation}
 is  split, up to multiplication by $[L_{\mathfrak{q}}:\Q_v]$, by corestriction,   and $[L_{\mathfrak{q}}:\Q_v]$ is invertible on $\mathcal{O}_{\mathfrak{p}}$.  
 
 This proves the first assertion, about $i$. For the second assertion, 
 note that the assumption about class groups means that
 $H^1_{\ur}(L, \mathbb{F}_{\mathfrak{p}}(1))$
  coincides with $\UnitL \otimes \mathbb{F}_{\mathfrak{p}}$.
The same analysis as above means that the rank of  $H^1_{\mathrm{ur}}(\Q, \Ad^* \rhobar(1))$ over
 $\mathbb{F}_{\mathfrak{p}}$ is bounded above by the dimension of
 $$ \left( \Ad^* \rhobar \otimes  \UnitL \right)^{G_{L/\Q}}$$
 and (again  because $G_{L/\Q}$ has no Galois cohomology in characteristic $p$)
 this dimension coincides with the $\mathcal{O}_{\mathfrak{p}}$-rank of $\OLrho \otimes \mathcal{O}_{\mathfrak{p}}$ (which is exactly $1$). 
The surjectivity now follows.  \qed

Now, under the identification $i:  H^1_{\mathrm{ur}}(\Q, \Ad^* \rho_{\mathfrak{p}}(1))  \simeq (\OLrho)_{\mathfrak{p}}$,
the composite
$$ H^1_{\mathrm{ur}}(\Q, \Ad^* \rho_{\mathfrak{p}}(1)) \rightarrow  H^1_{\mathrm{ur}}(\Q, \Ad^* \rhobar (1))  \stackrel{j^{\vee}}{  \rightarrow }H^1(\F_q, \Ad^* \rhobar(1))$$
$$
 \stackrel{\sim}{\rightarrow} H^1(\F_q, \mathbb{F}_{\mathfrak{p}}(1))  \simeq \F_q^*  \otimes (\mathcal{O}/\mathfrak{p})$$
is identified with the map $\theta_q$ described in \eqref{twmap0}.   Here we have made use of a map $\Ad^* \rhobar \rightarrow \mathbb{F}_{\mathfrak{p}}$, which comes from  pairing with the element defined in  \eqref{eqdef}. 
In particular, $\theta_q$ vanishes if and only if $j^{\vee}$ does, as required.

\section{Explication}  \label{sec:explicit}

Our main Conjecture     \ref{Mainconj},  as formulated, involves a cup product in coherent cohomology on the special fiber of a modular curve. We want to translate it to a readily computable form,
i.e. one that can be carried out just using manipulations with $q$-series. We will achieve this in this section, at least
in the case $n=1$ and under modest assumptions on $q$,  and then test the conjecture numerically.

\subsection{Pairing with $g'$}
Recall (\S \ref{gdef}) that we have fixed another weight one modular form $g'$ that is contragredient to $g$. 
To extract numbers from the Conjecture, we pair both sides of \eqref{equa} with $g'$, using the residue pairing (\S \ref{residuepairing}). 
Pairing \eqref{equa} with $g'$, and using 
$\theta_q^{\vee}(z) =  u^* \otimes  \frac{  \langle z, \theta_q(u) \rangle}{\langle u^*, u \rangle}$  from \eqref{theta vee explicit} 
we arrive at:

  \begin{equation} \label{doug2} [T_{q,z} \bar{g}, \bar{g}']_{\res, k} = \langle z, \theta_q(u) \rangle  \cdot   \overline{\left[ \frac{  \alpha [\theta_q^{\vee}(u) \star g, g']_{\res,\mathcal{O}}}{\langle u, u^* \rangle}\right]}\end{equation}
where both sides lie in $k$; and we recall again that we have written $\bar{g}$ for the reduction of $g$ to a modular form with $k$ coefficients. 

Now the square-bracketed quantity on the right hand side is an element of $E$, integral at $\mathfrak{p}$, and independent of choice of $(\mathfrak{p}, n, q, z)$. 
We abridge \eqref{doug2}
to 
 $$  [T_{q,z} \bar{g}, \bar{g}']_{\res, k} \propto   \langle z, \theta_q(u) \rangle $$
 This should hold true for any $(\mathfrak{p}, n, q, z)$. 
 
 Unwinding the definition of the derived Hecke operator, 
\begin{equation} \label{above2}  [T_{q,z} \bar{g}, \bar{g}']_{\res, k}  =  [\pi_{1}^* \bar{g} \cup z \Shim_X , \pi_2^* \overline{g'}]_{\res,k}\end{equation}
where the residue pairing is now taken on $X_{01}(qN)_k$,
 $\pi_1, \pi_2$ are the two projections $X_{01}(qN) \rightarrow X$,
and  $$z  \Shim_X  \in H^1(X_{01}(qN)_k,  \mathcal{O}  ).$$
(Recall that $\Shim_X \in H^1(X_{01}(qN)_k, \mathcal{O}\langle -1 \rangle)$, 
so its product with $z \in k \langle 1 \rangle$ lies in the right-hand group above.)
To simplify notation, define  the weight $2$ form 
\begin{equation} \label{GDef} G= \pi_1^* g  \cdot \pi_2^* g'   \in H^0(X_{01}(qN)_{k}, \Omega^1).\end{equation}
In terms of classical modular forms, $G$ would be the form ``$z \mapsto g(z) g'(qz)$.'' 
Then the right hand side of \eqref{above2} is simply the (Serre duality) pairing of $G \in H^0(\Omega^1)$
and $z \Shim_X \in H^1(\mathcal{O})$ in the coherent cohomology of $X_{01}(qN)$.  Therefore, the conjecture implies that
$\langle  z \Shim_X, G \rangle \propto \langle \theta_q(u), z \rangle$; and here we may as well cancel the $z$s from both sides:
\begin{equation} \label{doug3}  \langle  \Shim_X, G\rangle \propto  \theta_q(u). \end{equation}
 Here both sides lie in $k \langle 1\rangle$, that is to say, in $(\Z/q)^* \otimes k$.

Now the class $\Shim_X$ is pulled back from a class $\Shim$ on $X_0(q)$,
and correspondingly the pairing on the left-hand side can be pushed down to $X_0(q)$.
Writing 
$$  \Gproj = \mbox{ projection of $G$ to level $q \in H^0(X_0(q)_k,\Omega^1$)},$$
we have $\langle \Shim_X, G \rangle = \langle \Shim,  \Gproj \rangle$.

Thus our conjecture implies that  \begin{equation} \label{conj2} \langle \Shim,   \Gproj \rangle \propto  \theta_q(u), \mbox{ equality  in $(\Z/q)^* \otimes k$}. \end{equation}
where we recall that:
\begin{itemize}
\item $\Shim \in H^1(X_0(q)_{k},   \mathcal{O} \otimes (\Z/q)^* )$ is constructed from the covering $X_1(q) \rightarrow X_0(q)$;
\item $ \Gproj \in H^0(X_0(q)_{k} ,\Omega^1)$  is the pushforward of the form ``$z \mapsto g(z) g'(qz)$'' from level $X_{01}(qN)$  to level $X_0(q)$; it is a weight $2$ cusp form. 
\item $\langle -, -\rangle$ is the pairing of Serre duality. 
\item The symbol $\propto$ is interpreted as in \eqref{proptoequality}.
\end{itemize}

\subsection{Localization at the Eisenstein ideal}\label{Merel Eisenstein} 
 
 To translate \eqref{conj2} to a computable form, we will use computations of Merel and Mazur.
 Let  $$\mathrm{E} \in H^0(X_0(q)_{k} ,\Omega^1)$$ be the ``Eisenstein'' cusp form with $k$ coefficients, in other words,
the unique element  whose $q$-expansion coincides with the reduction modulo $p^n$ of the weight $2$ Eisenstein series;
the condition that $q \equiv 1$ modulo $p^n$ means that this weight $2$ Eisenstein series
indeed has cuspidal reduction in $k$. 
The pairing  
$$ \langle  \Shim, \mathrm{E} \rangle \in   (\Z/q)^*_p$$
 was considered by Mazur (\cite{mazur}, page 103, discussion of the element $u$)
 and was computed in a remarkable paper of Merel \cite{Merel}. 
 We will carefully translate
Merel's computation into our setting in the next section; unfortunately, in doing so,
we will have to impose the restriction $n=1$, i.e. we can only compute things modulo $p$
and not higher powers of $p$. 

\begin{lemma} \label{Merel lemma} (Merel; see \S \ref{sec:Merel} for details of the translation  
from Merel's framework to this one).  
\begin{equation} \label{Merelcomputation} \langle  \Shim, \mathrm{E} \rangle = \vM  \mbox{ mod $p$,}\end{equation}
where mod $p$ means that the two sides have the same projection to $\mathbb{F}_p\langle 1 \rangle$.\footnote{It seems likely that the two sides are actually equal in $(\Z/q)^*_p$ but we do not prove this.}

Here the {\em Merel unit} $\vM \in (\Z/q)^*$ is the element 
\begin{equation} \label{vMdef} \vM =  \zeta^2 \prod_{i=1}^{(q-1)/2} i^{-8i}, \ \ \  \zeta = \begin{cases} 1, q \equiv 2 (3), \\ 2^{(q-1)/3}, \mbox{ else} \end{cases}.\end{equation} 
\end{lemma}

In the remainder of this section, we will compute $\langle \Shim,   \Gproj \rangle$ 
(the left-hand side of \eqref{conj2})  using Lemma \ref{Merel lemma}.
 
 Let $\mathbb{T}$ be the Hecke algebra for cusp froms on $X_0(q)_{\Z_p}$, i.e.
the algebra of endomorphisms of
$S_2(q) := H^0(X_0(q)_{\Z_p}, \Omega^1)$
generated by $T_{\ell}$ for all $\ell \neq q$. 
 Let $\mathfrak{I} \leqslant \mathbb{T}$ be the Eisenstein ideal, i.e.
the kernel of the character
$$  \mathbb{T}  \rightarrow \Z/p\Z, \ \ T_{\ell} \mapsto (\ell +1) $$
by which $\mathbb{T}$ acts on the modulo $p$ reduction of $\mathrm{E}$. 
In particular it's a  maximal ideal.

Let $\mathfrak{m}_1, \dots, \mathfrak{m}_r$
be all other maximal ideals of $\mathbb{T}$. Then the natural map
$$ \mathbb{T} \longrightarrow \mathbb{T}_{\mathfrak{I}} \times \prod_{s=1}^{r} \mathbb{T}_{\mathfrak{m}_s}$$
is an isomorphism (here $\mathbb{T}_{\mathfrak{I}}$ means the completion, and similarly for $\mathfrak{m}_s$).   Let $e_{\mathfrak{I}}$ be the idempotent 
of $\mathbb{T}$ corresponding to the first factor; the splitting 
$$ 1 = e_{\mathfrak{I}}  +  \underbrace{ (1-e_{\mathfrak{I}})}_{=e_{\mathfrak{I}}'}$$
gives rise to a splitting
\begin{equation} \label{splitting} S_2(q) =  S_2(q)_{\mathfrak{I}} \oplus  S_2(q)_{\mathfrak{I}}'\end{equation}
where $S_2(q)_{\mathfrak{I}}$ is the  image of the idempotent $e_{\mathfrak{I}}$,
and the complementary subspace is the image of $1-e_{\mathfrak{I}}$. 
Therefore, if $T \in \mathfrak{I}$ is chosen so that $T \notin \bigcup_{i=1}^s \mathfrak{m}_i$,
then $T$ acts invertibly on the second factor.

  Decompose $\Gproj$ as
  $$ \Gproj = \Gproj_{\mathfrak{I}} + (\Gproj)'$$
  according to the splitting above. 
 The Shimura class $\Shim$  is annihilated by $\mathfrak{I}$ (see for example \cite{mazur}, Lemma 18.7).
 Choose as above $T \in \mathfrak{I}$ that acts invertibly on the second factor of \eqref{splitting}. We may write
 $$ \langle \Shim, (\Gproj)' \rangle = \langle \Shim, T T^{-1} (\Gproj)' \rangle = \langle T \Shim, T^{-1} (\Gproj)' \rangle = 0,$$ and so 
 $$ \langle \Shim,  \Gproj \rangle = \langle \Shim,  \Gproj_{\mathfrak{I}} \rangle$$
where as before the pairings come from Serre duality. 

Next, Mazur proves \cite[Proposition 19.2]{mazur}  that
\begin{equation} \label{Mazur nonvanish} \mbox{ $\vM$ is nonzero {\em modulo $p$}  $\iff$    
  $S_2(q)_{\mathfrak{I}}$ is of rank $1$ over $\Z_p$.} \end{equation}
  We will complete our computation only in this case.\footnote{Recently,  Lecouturier has proposed a very interesting generalization of  the conjectural equality \eqref{prediction}
  to the case when $\vM$ is zero modulo $p$ and has verified it numerically in some cases. }
  Since   $\mathrm{E}$ is annihilated by $\mathfrak{I}$, we have in fact
 $ \mathrm{E} \in S_2(q)_{\mathfrak{I}}$, and since the first Fourier coefficient of $\mathrm{E}$ is $1$, 
 we have  (under the assumption  of \eqref{Mazur nonvanish}) 
 $S_2(q)_{\mathfrak{I}} = \Z_p. E.$
Thus, after extending scalars to $\mathcal{O}$, we find 
$$   \Gproj_{\mathfrak{I}} = a_1( \Gproj_{\mathfrak{I}}) \cdot \mathrm{E},$$
 where $a_1( \Gproj_{\mathfrak{I}}) \in \mathcal{O} \otimes \Z_p$ denotes the first coefficient in the $q$-expansion.  
Putting this together  with our prior discussion, we have shown

\begin{Proposition} \label{prop:explicit}
Conjecture \ref{Mainconj}  implies that there exists $\alpha \in E$ such that 
\begin{equation} \label{prediction}  a_1( \Gproj_{\mathfrak{I}})  \otimes  \ (\vM)_p   \equiv \alpha \cdot  \theta_q(u) \mbox{ modulo  $p \cdot (\mathcal{O} \otimes (\Z/q)_p^*)$}\end{equation}
 for any $(\mathfrak{p}, n, q)$ as in \S \ref{lambdaconditions} with the additional property that
 $(\vM) \in (\Z/q)^*_p$ is nontrivial modulo $p$. \footnote{To be absolutely clear, we write out the meaning of this statement.   We understand $$ L := a_1( \Gproj_{\mathfrak{I}})  \otimes  \ (\vM)_p, R := \theta_q(u)$$
  as elements of $\mathcal{O} \otimes (\Z/q)^*_p$; and the statement above means
  that if we reduce $\bar{L}, \bar{R}$ to $\mathcal{O}/p \otimes (\Z/q)^*_p$, 
  then $\bar{L} = \alpha \bar{R}$, in the sense of \eqref{prop test}. }
  Other conventions are as follows:

 \begin{itemize}
 \item $\vM \in (\Z/q)^*_p$ is the Merel unit, see \eqref{Merelcomputation}; 
 \item $a_1( \Gproj_{\mathfrak{I}}) \in \mathcal{O} \otimes \Z_p$ is the first Fourier coefficient of $G = (\pi_1^* g) (\pi_2^* g')$, after taking projection $ \Gproj$ to level $X_0(q)$
 and then projection $\Gproj_{\mathfrak{I}}$ to the localization at the Eisenstein ideal.   
 \item $\theta_q(u) \in  k \langle 1 \rangle  = \mathcal{O} \otimes (\Z/q)^*_p$ is the reduction of the Stark unit
 \end{itemize}
 
  \end{Proposition}

 \subsection{Some philosophical worries}
 
 Let us take  to examine some consequences of an inadequacy of our conjecture, namely, 
  it is only  formulated ``up to $E^*$'.''
   
For each $(\mathfrak{p}, n, q)$ as in \S \ref{lambdaconditions}, 
we can compute both $ a_1( \Gproj_{\mathfrak{I}})  \otimes  \  (\vM)_p$
and $ \theta_q(u)$ and compare them.  Let us also restrict 
 to  $(\mathfrak{p}, n, q)$ for which $\theta_q(u) \neq 0$; there are infinitely many such $\mathfrak{p}$. 
Therefore \eqref{prediction} specifies the reduction $\alpha \in E$ to $\mathbb{P}^1(\mathbf{F}_{\mathfrak{p}})$,
for an infinite collection of $\mathfrak{p}$. This uniquely specifies $\alpha$ if it exists.
  
 The conjecture is numerically falsifiable to some extent.
 For example, if we find two different pairs $(\mathfrak{p},n, q)$
 and $(\mathfrak{p}, n', q')$ 
 for which the predicted reductions of $\alpha$ mod $\mathfrak{p}$ differ,
 this clearly contradicts the conjecture.  
Indeed the fact that this did not occur in our numerical computations was very encouraging to us. 
 
 However,
if this type of clash does not occur, no amount of computation can falsify the conjecture:
we can, of course, produce an $\alpha \in E$ with any specified reduction at any number of places. 
Nonetheless this proves to be  largely a theoretical worry. 
 In our examples, we shall find an $\alpha$ of very low height
 for which \eqref{prediction} holds for many $(\mathfrak{p}, n, q, z)$. 
 Our sense is that this should be taken as    satisfactory indication
 that the Conjecture, or something very close to it at least, is valid. 
 
As a final excuse we may note that the  conjectures about special values of $L$-functions
 were initially phrased  with a $\Q^*$ ambiguity that  is similarly unfalsifiable. 

Eventually, we hope that these issues will be solved by formulating an integral
  form of the conjecture; this could perhaps be done using the theory
  of derived deformation rings.

  \subsection{Forms associated to cubic fields}  \label{sec:cubic}
 
We now make the foregoing discussion even more explicit for  the form $g$ associated to a cubic field $K$;
write $L$ for the Galois closure of $K$. (This will coincide with our previously defined $L$ in a moment.)

 Such a field $K$ defines a representation $\Gal(L/\Q) \rightarrow S_3$;
if we regard $S_3$ as acting on  $M =  \{(x_1, x_2, x_3) \in \Z^3: \sum x_i = 0\}$ by permuting the coordinate axes, 
we may regard $\rho$ as a rank $2$ Galois representation:
\begin{equation} \label{browncow}  \rho: \Gal_{\Q} \twoheadrightarrow S_3 \rightarrow \GL_2(M). \end{equation}

  Under the representation \eqref{browncow}, there is a basis for $M$ such that  the transposition $\sigma = (12) \in S_3$ is sent to  $S:=\left( \begin{array}{cc}  0 & 1 \\ 1 & 0  \end{array}\right)$,
whereas a $3$-cycle  $\tau= (123) \in S_3$ is sent to $T := \left( \begin{array}{cc}   -1 & 1 \\ -1 & 0  \end{array}\right)$.  
We may set things up so that the fixed field of $(12) \in S_3$ is equal to $K$.

In our previous notation, take 
\begin{itemize}
\item $L$ as above, namely, the Galois closure of the cubic field $K$;
\item  $E=\Q$ and $\mathcal{O} = \Z$. 
\item $\mathfrak{p}  = p \geqslant 5$ to be a rational prime of $\Q$.
\item $n=1$ (thus we work only modulo $p$ rather than $p^n$). \item $q \equiv 1 (p)$ to be a prime such that the $q$th Hecke eigenvalue $a_q(g) = 0$. 
In this case,  the Frobenius is a transposition\footnote{The primes $q$ for which $\rho(\Frob_q)$ is a $3$-cycle also are Taylor--Wiles primes, but it is then easy to see that $T_{q,z} g = 0$ for such $q$.
To verify this, one can use the fact -- notation as in \eqref{GDef} --  that the Atkin-Lehner involution at $q$ for $X_{01}(qN)$
acts by $-1$ on $\Shim_X$, but it acts by $\chi(q)$ on $G$, where $\chi$ is the quadratic Nebentypus character for $g$. }  
 in $S_3$. Thus  $q$ 
is a Taylor--Wiles prime with eigenvalues $(1,-1)$. 
\item Therefore in this case
$k \langle 1 \rangle = \mathbb{F}_p\langle 1 \rangle$  is just the unique quotient of $(\Z/q)^*$ of order $p$. 

\item We also 
fix a prime $\mathfrak{q}_0$ of $L$ over $q$ such that the image of the Frobenius for $\mathfrak{q}_0$ is equal to $S$. 
In particular, this fixes $K$, so the prime $\tilde{q}$ of $K$ below $\mathfrak{q}_0$ is of degree $1$ over $q$. 
 \end{itemize}

\begin{lemma}
Consider the isomorphism $\OLrho \simeq \Hom_{G_{L/\Q}}(\Ad^0 \rho, \UnitL)$ of  \eqref{dododo}. 
(Recall that $\UnitL$ is the unit group of $L$.)
Then computing the image of $S \in \Ad^0 \rho$  
gives rise to an isomorphism  \begin{equation}
\label{OLexplicit} \OLrho  \otimes \Z[\frac{1}{6}] \simeq \mathcal{O}_K^{(1)} \otimes \Z[\frac{1}{6}],\end{equation}
where $\mathcal{O}_K^{(1)}$ is the group of norm one units of $K$. 

Moreover, 
 for $p \geq 5$ 
 the reduction map $\theta_q: \OLrho \rightarrow\mathbb{F}_p \langle 1 \rangle$ described
in \eqref{twmap0}  becomes identified with the composite
$$  \mathcal{O}_K^*  \rightarrow     (\mathcal{O}_K/\tilde{q})^* = (\Z/q)^* \rightarrow \mathbb{F}_p \langle 1 \rangle,$$
where $\tilde{q}$ is the unique degree one prime of $K$ above $p$. 
\end{lemma}
\proof 
Indeed we may split
$$\Ad^0 \  \rho  \otimes \Z[\frac{1}{6}] =  \Hom^0(M,M) \otimes \Z[\frac{1}{6}] = \Z[\frac{1}{6}]   e   \oplus \mathsf{W} $$
where $e$ is the projection of $T \in \Hom(M,M)$ to the trace zero subspace $\Hom^0$,  and $\mathsf{W}$ is 
the $\Z[\frac{1}{6}]$-submodule of $\Hom(M, M) \otimes \Z[\frac{1}{6}]$  spanned  by the images of $(12), (13), (23)$ under $\rho$. 

Therefore $S_3$ acts on $e$ by the sign character,
whereas  for any $S_3$-module $V$, the space of homomorphisms $\Hom_{S_3}(\mathsf{W}, V)$ is identified with the subspace of $$v \in V^{(12)}=\mbox{$(12)$-fixed vectors in $V$}$$ such that $v+(123)v+ (132) v=0$.

    Using the definition of $\OLrho$
and the splitting above, we find that evaluation at $S$ induces an isomorphism
$$ \OLrho \otimes \Z[\frac{1}{6}] \simeq \left(  \UnitL^{(\mathrm{sign})} \oplus \mathcal{O}_K^{(1)}\right) \otimes \Z[\frac{1}{6}].$$
The first factor corresponds to units in the imaginary  quadratic field $\Q(\sqrt{\disc(L)})$, and 
is thus trivial upon inverting $6$.   This proves \eqref{OLexplicit}.
  
Now let $u$ be a norm one unit in $K$; we may now identify it with an element of $\OLrho \otimes \Z[1/6]$. We will compute its image under the reduction map.   
Let  $\mathbf{u} \in \Hom(\Ad^0 \rhobar,  \UnitL )$ be the element associated to $u$.
By definition  $\mathbf{u}(S) =u$. Let $\mathfrak{q}_0$ be the prime of $L$ above $\tilde{q}$, as before;  
to compute $\theta_q(u)$ we must, by definition, 
  compute the image of $\mathbf{u}$ under the sequence \eqref{twmap0}:
$$ \Hom_{G_{L/\Q}}(\Ad^0 \rhobar,  \prod_{\mathfrak{q}|q} \mathbf{F}_{\mathfrak{q}}^* ) \stackrel{\sim}{\rightarrow}  \Hom(\Ad^0 \rhobar,  \mathbf{F}_{\mathfrak{q}_0}^*)^{D_{\mathfrak{q}_0}} \stackrel{\mathsf{e}_q}{\longrightarrow}  \mathbb{F}_p \langle 1 \rangle$$
where we phrased the previous definition dually.  The element $\mathsf{e}_q$ from \eqref{eqdef} is identified here with $S$, so that the last map is evaluation at $S$. 
It follows that this map is simply the reduction of $u$ at $\tilde{q}$.

 \qed

It follows from this discussion and Proposition \ref{prop:explicit} that we can rephrase our conjecture
 in the following way:
 
 \begin{Conjecture} Let $K$ be a cubic extension with  negative discriminant $-D$,
 with sextic Galois closure $L$. 
Let $g$ be the associated weight one form of level $D$.  
 Let $u \in \mathcal{O}_K^*$ be a unit.   Let $q \equiv 1$ modulo $p$ be as above; suppose that
    $\left( \frac{-D}{q} \right) = -1$, and $p \geq 5$, 
   and finally $\vM \in (\Z/q)^*$ (see \eqref{vMdef} for definition) is nonzero 
   modulo $p$, i.e. upon projection to the quotient $\mathbb{F}_p \langle 1 \rangle$.

Then  there exist  $A, B \in \mathbf{Z}$ such that, for all such $q$
 we have   \begin{equation}   \vM^ {A \cdot   \eta   }  =      \bar{u}^B \mbox{  in  $ \mathbb{F}_p \langle 1 \rangle$.} \end{equation}
  where:
  \begin{itemize}
  \item $\eta \in \Z $ is the   first Fourier coefficient of the Eisenstein component of  $ \Gproj_{\mathfrak{I}}$, the projection of $g(z) g(qz)$
  to the Eisenstein component at level $q$. (This is well defined modulo  the numerator of $\frac{q-1}{12}$, which is sufficient to make sense of the  above definition.) 
\item  $\bar{u} \in (\Z/q)^*$ is the reduction of $u$ modulo the unique degree one prime of $K$, above $q$. 
\end{itemize}

 \end{Conjecture}  

 We have tested this conjecture numerically (see data tables) for the fields $K$ of discriminant $-23$ and $-31$.
In all the cases for discriminant  $-23$ we find $\frac{A}{B} = \frac{-1}{72}$; in all the cases for discriminant $-31$
we find $\frac{A}{B} = 72$. The fact that $72$ is divisible only by $2$ and $3$ is striking.

 \begin{table}[!htbp]
  \begin{tabular}{|c|c|c|c|c|}
 \hline \hline
 $p$ &   $q$   & $\log(\bar{u})/\log(\vM) \in \Z/p$ & $\eta \in \Z/p$  & ratio \\ 
 \hline
 $5$ & $11$   & 3 (5)  &  4(5)  & 2(5) \\ 
 $5$ & $61$    &  1(5)  & 3(5 & 2(5)) \\ 
 $5$ & $81 $ & $\infty$ & - & - \\ 
 $7$ & $43$ & 3(7) &  1(7)  & 3(7) \\
  $7$ & $113$ & 1(7)   & 5(7)   & 3(7)   \\
 $11$ & $67$        & 6(11)  & 8 (11)  & -2(11)\\ 
 $11$ & $89$       & 1(11) & 5(11)  & -2(11)  \\ 
 $13$ & $53$     &  6  (13) &   10(13) & -2(13) \\ 
 $13$ & $79 $     & 5 (13) &   4 (13)  & -2 (13)\\
  $17$ & $103 $     &   $\infty$ & - & -    \\
  $17$ & $137$ &  5(17) &14 (17) & 4 (17)  \\
  $37$ & $149$ &  20(37) & 3 (37) & 19 (37)  \\
  $41$ & $83 $     & 12(41) &   38(41) & -4(41) \\
    $53$ & $107 $     &   30(53) &  13(53)   & 39 (53) \\
 \hline \hline 
 \end{tabular}
 \caption{Data for the weight one form associated to the cubic field with discriminant $-23$; in all cases the ratio is $-1/72$ modulo $p$.
 All allowable $p\leq 100$ and $q \leq 150$ shown. } \label{Table1}
   \end{table}

  \begin{table}[!htbp]
  \begin{tabular}{|c|c|c|c|c|}
 \hline \hline
 $p$ &   $q$   & $\log(\bar{u}/\vM) \in \Z/p$ & $\eta \in \Z/p$  & ratio \\ 
 \hline
 $5$ & $11$   & 2 (5)  &  4(5)  & 2(5) \\ 
 $5$ & $61$    &  2(5)  & 4(5 & 2(5)) \\ 
  $7$ & $29$    &  1(7)  & 2(7) & 2(7)) \\ 
 $7$ & $43$    &  4(7)  & 1(7) & 2(7)) \\ 
  $7$ & $127$    &   $\infty$  & - & - \\ 
$11$ & $23$ & 3 (11) & 7 (11) & 6(11)  \\
$11$ & $89$ & 7(11) & 9(11) & 6(11) \\
$13$ & $53$  & 2(13) & 1(13) & 7(13) \\
$13$ & $79$ & 3(13) & 8 (13)&  7 (13) \\
$17$ & $137$ & 4(17) & 16 (17)  & 4(17) \\ 
$23$ & $139$  & 4(23) & 12 (23) & 3 (23) \\
 $41$ & $83$ &  28(41) &7(41) &31(41)\\
  \hline \hline 
 \end{tabular}
 \caption{Data for the weight one form associated to the cubic field with discriminant $-31$; in all cases the ratio is $72$ modulo $p$. 
 All allowable $p \leq 100$ and $q \leq 150$ shown. ? means that we did not compute because it took too long; - means
 undefined.  } \label{Table2}
   \end{table}

   \section{Flat cohomology and Merel's computation} \label{sec:Merel}

We now explain why Merel's computation implies Lemma \ref{Merel lemma}. The issue is that   
  Merel's computation is in characteristic zero. To relate it to $\langle \mathrm{E}, \Shim \rangle$, which is defined
  in characteristic $p$, we will need to do a little setup in flat cohomology.

  Let $X=X_0(q)$ regarded now as a a proper smooth curve over $\Z_p$;
  here $q \equiv 1$ modulo $p$.  Let $J_p$ be the $p$-torsion of the Jacobian of $X_{\overline{\Q_p}}$. 
We shall define several incarnations of both the Shimura class and the Eisenstein class. 

\subsubsection{The (Shimura) class $\alpha$}

 The Shimura cover $X_1(q)^{\Delta} \rightarrow X_0(q)$ (from \S \ref{Shimuraclassdef}) is a $(\Z/q)^*_p$ torsor for the etale topology.
As before it define a class  $\Shim \in H^1_{\et}(X, \mathbb{F}_p \langle 1 \rangle)$,
which can be pulled back to flat cohomology:
$$ \alpha \in H^1_{\fl}(X,  \mathbb{F}_p \langle 1 \rangle).$$
Restricting $\Shim$ to the geometric generic fiber $X_{\overline{\Q_p}}$  we get a class in {\'e}tale cohomology
$$\alpha_{\et} \in H^1_{\et}(X_{\overline{\Q_p}}, \mathbb{F}_p \langle 1 \rangle).$$

  The inclusion $\mu_p \hookrightarrow \mathbb{G}_m$ induces
 $ H^1_{\et}(X_{\overline{\Q_p}}, \mu_p) \rightarrow   J_p$,
and thus $\alpha_{\et}$ gives
\begin{equation} \label{Palphadef} P_{\alpha} \in \Hom( \mu_p \langle -1 \rangle, J_p)\end{equation}
(we use the notation $P_{\alpha}$ to suggest that this is a point on the Jacobian). 

Finally, we also obtain a Zariski class on the geometric special fiber, using
the inclusion $\mathbb{F}_p \hookrightarrow \mathcal{O}$ and the identification of Zariski and etale cohomology for $\mathcal{O}$: 
$$ \alpha_{\Zar} \in H^1_{\Zar}(X_{\overline{\mathbb{F}_p}}, \mathcal{O} \langle 1 \rangle).$$

\subsubsection{The (Eisenstein) class $\beta$}

Let $\Delta$ be the weight  $12$ cusp form $q \prod (1-q^n)^{24}$, and consider the function $f := \Delta(qz)/\Delta(z)$ on $X$.
  Extracting its $p$th root gives a $\mu_{p}$-torsor (in the flat topology) on $X$. Indeed, 
 $f$ is invertible except for the divisors corresponding to $0$ and $\infty$, 
 and along those divisors its valuation is divisible by $p$.   Thus, we get a class
 $$\beta \in H^1_{\fl}(X, \mu_{p}).$$
 
 The $\mu_{p}$-torsor  is {\'e}tale over
the geometric generic fiber $X_{\overline{\Q_p}}$ and we 
we get a corresponding class in {\'e}tale cohomology 
$$\beta_{\et} \in H^1_{\et}(X_{\overline{\Q_p}}, \mathbb{F}_p \langle 1 \rangle).$$
  
 There is a corresponding class in the $p$-torsion of the Jacobian, namely,
 writing $0$ and $\infty$ for the two cusps of $X$ we may form
 $$ Q_{\beta} :=  \frac{(q-1)}{p} ( (\infty) - (0)) \in J_p$$
-- this is related to our prior discussion because $p. Q_{\beta}$ is the divisor of $f$. 
  
Finally, there is also a Zariski class ``corresponding'' to $\beta$ on the special fiber. Namely, 
  the logarithmic derivative $\frac{df}{f}$ in fact extends to a global section of $\Omega^1$, i.e. a class
$$\beta_{\Zar} \in H^0(X_{\overline{\mathbb{F}_p}}, \Omega^1).$$
Observe that   $\frac{df}{f}$ is the differential form associated to the  ``Eisenstein cusp form'' $G$ of weight $2$.

With these preliminaries,  the main point is to check the following

\begin{Proposition} We have an equality in $\mathbb{F}_p\langle 1 \rangle$: 
$$\langle P_{\alpha}, Q_{\beta} \rangle_{\mathrm{Weil}} =  \langle \alpha_{\et}, \beta_{\et} \rangle_{\et} = \langle \alpha_{\Zar}, \beta_{\Zar} \rangle_{\Zar}.$$
\end{Proposition}
Here $\langle -, - \rangle_{\mathrm{Weil}}$ is the Weil pairing, $\langle - , - \rangle_{\et}$ is the pairing given by Poincar{\'e} duality in {\'e}tale cohomology on the geometric fiber, and 
$ \langle - , - \rangle_{\Zar}$ is the pairing given by Serre duality in coherent cohomology on the special fiber.   Keeping track of twists
we see that these all take values in $\mathbb{F}_p \langle 1 \rangle$. 

 Now  $ \langle \mathrm{E} , \Shim \rangle$ is given by $\langle \alpha_{\Zar}, \beta_{\Zar} \rangle_{\Zar}$; the Proposition shows this coincides  
 (in $\mathbb{F}_p \langle 1\rangle$) with $ \langle  P_{\alpha}, Q_{\beta} \rangle_{\mathrm{Weil}}$. The  Weil pairing on the right is computed by Merel; 
we pin down the relation to Merel's computation in \S \ref{Weil Merel}.  Taken together, the Proposition 
and this computation prove Lemma \ref{Merel lemma}.
 
 \proof  
 The first equality is straightforward:  an explicit representative for $Q_{\beta} \in J_p \simeq H^1(X, \mu_p)$ 
is given by the $\mu_p$-torsor  associated to $f :=  \Delta(qz)/\Delta(z)$, 
because the divisor of $f$ is $p Q_{\beta}$.

 We now discuss the second equality.   We will compare everything to    the cup product in flat cohomology, i.e.  $$ \alpha \cup \beta \in H^2_{\fl}(X, \mu_{p} \langle 1 \rangle).$$

  There is a degree map $H^2_{\fl}(X, \mu_p) \rightarrow \mathbb{F}_p$; let us explicate it. 
On any scheme, the sequence  
  $\mu_{p} \rightarrow \mathbb{G}_m \rightarrow \mathbb{G}_m$
induces an exact sequence of represented sheaves for the flat topology.
This identifies the flat cohomology of $\mu_{p}$
with the hypercohomology of $[\mathbb{G}_m \stackrel{x \mapsto x^p}{\longrightarrow} \mathbb{G}_m]$. 

Let $X_{\overline{\Z_p}}$ be the base change of $X$ to $\overline{\Z_p}$
 (the Witt vectors of $\overline{\mathbb{F}_p}$).  We obtain an exact sequence
 \begin{equation} \label{ex21}  \mathrm{Pic}(X_{\overline{\Z_p}})/p \hookrightarrow  H^2_{\fl}(X_{\overline{\Z_p}}, \mu_{p}) \rightarrow H^2_{\fl}(X_{\overline{\Z_p}}, \mathbb{G}_m)[p]\end{equation}
Flat and \'etale cohomology of $\mathbb{G}_m$ coincide (see \cite{***}), and
 the right-hand side is a subgroup of the Brauer group of $X_{\overline{\Z_p}}$,
which  vanishes  (Th{\'e}or{\`e}me 3.1 of \cite{groth}).  
 Accordingly, any class in $H^2_{\fl}(X_{\overline{\Z_p}}, \mu_p)$
 is the coboundary of a line bundle, and computing degree gives the desired homomorphism
 $$ \mathrm{deg}: H^2_{\fl}(X_{\overline{\Z_p}}, \mu_p) \longrightarrow \mathbb{F}_p.$$
 We see that 
 $ \mathrm{deg} (\alpha \cup \beta) = \langle \alpha_{\et}, \beta_{\et} \rangle_{\et}$  and so it remains to see 
 $$ \mathrm{deg}(\alpha \cup \beta) = \langle \alpha_{\Zar}, \beta_{\Zar} \rangle_{\Zar}.$$

  Let $\pi$ be the morphism from the flat site on $X_{\overline{\mathbb{F}_p}}$ to the {\'e}tale site.
  As a reference for what follows, we refer to the paper of Artin and Milne \cite{ArtinMilne}.     We have isomorphisms:
 $$ R\pi_* (\Z/p\Z) \simeq [\mathcal{O} \stackrel{1-\mathrm{F}}{\rightarrow} \mathcal{O}]$$
 $$ R\pi_* \mu_p \simeq [\Omega^1 \stackrel{1-\mathrm{C}}{\rightarrow} \Omega^1][1]$$
 where $\mathrm{F}$ and $\mathrm{C}$ are, respectively, the Frobenius and Cartier maps. 
 and Artin--Milne show that the pairing $\Z/p\Z \times \mu_p \rightarrow \mu_p$
 induces, after push-forward, the ``obvious'' pairing on the complexes on the right, which can be computed in the Zariski topology, because flat and Zariski cohomology coincide for quasi-coherent sheaves.  
 
For the same reason,  the second identification induces an isomorphism 
 $$ H^2_{\fl}(X_{\overline{\mathbb{F}}_p}) \simeq H^1(X_{\overline{\mathbb{F}}_p}, \Omega^1)^{\mathrm{C}} = \mathbb{F}_p,$$
 where the map $H^1(\Omega^1) \rightarrow \overline{\mathbb{F}_p}$ comes from Serre duality.  Moreover, the resulting identification is simply the degree map, alluded to above; 
 this comes down to the fact that the map
 $$H^1(X_{\overline{\mathbb{F}}_p}, \mathbb{G}_m) \stackrel{d\log}{\longrightarrow} H^1(X_{\overline{\mathbb{F}}_p}, \Omega^1)  \rightarrow \overline{\mathbb{F}_p}$$
 again computes the degree of a line bundle modulo $p$.

With respect to the resulting identification
of $H^1_{\fl}(X, \mu_p) \simeq \mathbb{H}^0(\Omega^1 \stackrel{1-\mathrm{C}}{\rightarrow} \Omega^1)$,
and the \v{C}ech representation of this last hypercohomology, the class $\beta$ is represented by $\frac{df}{f}  \in \check{C}^0(\Omega^1)$, which has zero boundary and which is annihilated on the nose by $1-\mathrm{C}$.
 Similarly the class $\alpha$ in {\'e}tale cohomology is represented by a Cech cocycle $c^1 \in \check{C}^1(\mathcal{O})$ together with 
 a class $c^0 \in \check{C}^0(\mathcal{O})$ satisfying
 $ (1-\mathrm{F})c^1 = d c^0$. The image of the pairing  $\alpha \cup \beta \in H^2_{\fl}(\mu_p)$,
 under the map $H^2_{\fl}(\mu_p) \rightarrow H^1(\Omega^1)^{\mathrm{C}}$, is represented by 
 $c^1 \cdot \frac{df}{f} \in \check{C}^1(\Omega^1)$; its image by the trace pairing is the usual Serre duality pairing between the cohomology classes of $c^1$ and $\frac{df}{f}$.
  This concludes the proof.
  \qed

\subsection{Merel's computation} \label{Weil Merel}
Although routine, we write out  the details involving $\langle P_{\alpha} , Q_{\beta}\rangle$  to be sure of factors involving $\gcd(q-1,12)$. 
In what follows, we understand our modular curves to be considered  over an algebraically closed field of characteristic zero. 

 Recall that $P_{\alpha}$ is an element of
 $ \Hom( \mu_p \langle -1 \rangle, J_p)$. 
Thus the Weil pairing $ \langle P_{\alpha} , Q_{\beta} \rangle \in  \mathbb{F}_p \langle 1 \rangle$
  has the property that  \begin{equation} \label{previous}   \mbox{Weil pairing of } P_{\alpha}(u)  \mbox{ and }Q_{\beta}  =  u \cdot  \langle P_{\alpha}, Q_{\beta} \rangle_{\mathrm{Weil}} \ \ \ (u \in   
\mu_p\langle -1 \rangle). \end{equation}
where, on the  left hand side we have the ``usual'' Weil pairing of two torsion points in $J_p$.

  Following Merel, let $\nu$ be the gcd of $q-1$ and $12$; let $n = \frac{q-1}{\nu}$. 
Let $U \subset (\Z/q)^*$
be the subgroup of $\nu$th powers; the map
  $(\Z/q)^* \rightarrow \mathbb{F}_p \langle 1 \rangle$ factors through the $\nu$th power map, and we get a sequence
   $$ (\Z/q)^* \stackrel{x \mapsto x^{\nu}}{ \longrightarrow} U  \rightarrow \mathbb{F}_p \langle 1 \rangle.$$

The Galois group of the covering $X_1(q) \rightarrow X_0(q)$  can be identified
with $U$ (as in \S 3.3 \cite{Merel}).  This gives rise to a map
$$  \alpha': \Hom(U, \mu_{n}) \rightarrow  J_n$$
Also $Q' = (\infty) - (0)$ is $n$-torsion in the divisor class group, thus defining another class in $J_n$. 
Then Merel shows that 
   \begin{equation} \label{merelshows} \langle \alpha'(t), Q' \rangle_{n} = t (\vM), \ \ t \in \Hom(U,\mu_n)\end{equation}
    where the equality is in $\mu_n$ and 
  the subscript $n$ means we are using the Weil pairing at the $n$-torsion level. 

   We want to compare $\alpha'$ to $P_{\alpha}$.  
   Note that if $t \in \Hom(U, \mu_n)$ the power $t^{n/p}$ 
   defines an element of $\Hom(U, \mu_p)$ which, considered as an element
   of $\Hom((\Z/q)^*, \mu_p)$, factors through $\mathbb{F}_p \langle 1 \rangle$.  We refer to the resulting element
   as $\bar{t} \in \Hom(\mathbb{F}_p \langle 1 \rangle, \mu_p)$.   Explicitly, if $\mu \in (\Z/q)^*$, we have 
\begin{equation} \label{confusing} t^{n/p}(\mu^{\nu}) = \bar{t}(\mu)\end{equation}

Now  consider the commutative diagram (where we write $X=X_0(q)$ for short)   
  \begin{equation*}
    \xymatrix{
   H^1(X, U ) \times  \Hom(U, \mu_n) \ar[r] \ar[d]^{\mathrm{id} \times t \mapsto \bar{t} }&  J_n \ar[d]^{\times n/p}\\
     {} \underbrace{ H^1(X, U)  }_{\rightarrow H^1(X,\mathbb{F}_p\langle 1\rangle)}\times   \underbrace{ \Hom(\mathbb{F}_p\langle 1\rangle, \mu_p)}_{\simeq \mu_p \langle - 1\rangle } \ar[r]  &  J_p }
  \end{equation*}
  When we evaluate at the element of $H^1(X,U)$ corresponding to the cover $X_1(q) \rightarrow X_0(q)$, the
  top horizontal map becomes $\alpha'$ and the bottom map becomes $P_{\alpha}$ from
  \eqref{Palphadef}. 
Thus we have
$$\alpha'(t)^{n/p} = P_{\alpha} (\bar{t}), \ \ t \in \Hom(U, \mu_n) \mapsto \bar{t} \in \Hom(\mathbb{F}_p \langle 1 \rangle, \mu_p)$$ 
Pairing with $Q_{\beta} = \frac{q-1}{p} Q' \in J_p$  and comparing with \eqref{previous}:

$$ \underbrace{ \bar{t}}_{\mu_p(-1)}  \underbrace{ \langle P_{\alpha}, Q_{\beta} \rangle}_{\mathbb{F}_p \langle 1 \rangle} = \langle P_{\alpha}(\bar{t}), Q_{\beta} \rangle_p  = \langle \alpha'(t)^{n/p}, \frac{q-1}{p} Q' \rangle_p = \frac{q-1}{p} \langle \alpha'(t), Q' \rangle_n \in \mu_p$$
and so
$$ \underbrace{ \bar{t}}_{\mu_p(-1)}  \underbrace{ \langle P_{\alpha}, Q_{\beta} \rangle}_{\mathbb{F}_p \langle 1 \rangle} \stackrel{\eqref{merelshows}}{=}     \frac{q-1}{p} t(\vM) 
=    t^{n/p}(\vM^{\nu})\stackrel{\eqref{confusing}}{=}     \bar{t}(\vM),$$ 
where the equality is once again in $\mu_p$. 
We conclude that  $\langle P_{\alpha}, Q_{\beta} \rangle $ is indeed the  image of $\vM$ inside $\mathbb{F}_p \langle 1 \rangle$.

   \section{Comparison with the theory of \cite{V3}} \label{final}
  
 Derived Hecke operators at Taylor-Wiles primes have been defined abstractly  for general $q$-adic groups in \cite{V3}.
  The purpose of the present section is to identify the operators introduced in \ref{dhdefn} with those defined in \cite{V3}.    
  (The results of this section are, strictly speaking, not used elsewhere in the paper; however they show that
  all the constructions we have made are inevitable.)

Write $G = \GL_2(\Q_q)$ where $q \equiv 1 \pmod{p}$,
and $K = \GL_2(\Z_q)$. 
  Fix a base ring $S$ that is a $\Z_p$-algebra.

What we will need to do, in order to study   the derived Hecke operator at $q$, is to identify the cohomology of the modular curve with  the cohomology of the  $K$-invariants of a complex
of $G$-representations.  Unsurprisingly,  this is done by adding infinite level at $q$; we just pin down the details. 
We need to take a little care because the tower of coverings that one gets by adding infinite $q$-level is not {\'e}tale; however,
its ramification is prime to $p$, which will be enough for our purposes. 

In particular, we will use\footnote{with an apology to 21st century readers, see below...}  Lemma A.10 of  Appendix A of \cite{V3},
which  explicates the    action of the abstract derived Hecke algebra in terms of restrictions, corestrictions, and cup products. 

\subsection{Construction of complexes with an action of $\GL_2(\Q_q)$.}

 Let us fix a level structure away from $q$ for the usual modular curve, i.e., an open compact subgroup
 $K^{(q)} \subset \GL_2(\mathbb{A}^{(\infty, q)})$. We require that $K^{(q)} = \prod_{v \neq q} K_v$, where $K_v$ is hyperspecial maximal
 for almost all $v$.

For $U \subset \GL_2(\Q_q)$ an open compact subgroup, let $X(U)$ 
be the Deligne--Rapoport compactification of the modular curve with level structure $K^{(q)} \times U$.
 This again   has (Deligne--Rapoport) a smooth proper model over $\Spec \ S$, denoted $X(U)_S$.
 We denote again by $\CL_U \rightarrow X(U)$ the
 relative cotangent bundle of the universal elliptic curve; this defines a locally free sheaf over $X(U)_S$.

 Let us consider the pro-system of schemes $$X_{\infty} \ : U \mapsto X(U)$$
 indexed by the collection of all open compact subgroups of $\GL_2(\Q_q)$; the maps are inclusions $V \subset U$ of open compact subgroups.

 The isomorphisms
$X(g^{-1} U g) \stackrel{\sim}{\rightarrow} X(U)$ induce an action of $G = \GL_2(\Q_q)$ on $X_{\infty}$ (considered
as a pro-object in the category of schemes).  Let $\CL_{\infty}$ be the ``vector bundle'' over $X_{\infty}$ defined by $\CL$: by this we mean that   $\CL_{\infty}$ 
is a pro-scheme over $X_{\infty}$, which is level-wise a vector bundle.
 
 We will need the following properties:
\begin{itemize}
\item[(i)]  The action of $G$ on $X_{\infty}$ lifts to an action on $\CL_{\infty}$.  
\item[(ii)] Suppose that $V$ is a normal subgroup of $U$. Then
the natural map $$f_{UV}: X(V)_S \rightarrow X(U)_S$$
is finite, and identifies $X(U)_S$ with the quotient of $X(V)_S$ by $U/V$ in the category of schemes.    (See \cite[3.10]{DR}).  

Moreover,   there is a natural (in $S$) identification 
  $f_{UV}^* \omega_U \simeq \omega_V$.

\item[(iii)] With notation as in (ii), if the order of $U/V$ is a power of $p$, then the map $X(V)_S \rightarrow X(U)_S$ is {\'e}tale.  
 \end{itemize}

\proof (of (iii) only:) We may suppose that $S=\Z_p$.   The map is {\'e}tale over the interior of the modular curve,
so, by purity of the branch locus, it is enough to check that it is {\'e}tale at the cusps in characteristic zero.
The cusps of a modular curve are parameterized by an adelic quotient, but replacing
the role of an upper half-plane by $\mathbb{P}^1(\Q)$; so 
 we must verify that the map $$ \GL_2(\Q) \backslash (\GL_2(\mathbb{A}_f) \times \mathbb{P}^1(\Q))/V \longrightarrow  \GL_2(\Q) \backslash (\GL_2(\mathbb{A}_f) \times \mathbb{P}^1(\Q))/U,$$
considered as a morphism of groupoids, induces isomorphisms on each isotropy group.

Let $\mathbf{B}$ be a Borel subgroup in $\GL_{2,/\Q}$ and $\mathbf{N}$ its unipotent radical. 
We can identify $\mathbb{P}^1(\Q)$ with $\GL_2(\Q)/\mathbf{B}(\Q)$. 
The desired result follows, then, if 
for each $g \in \GL_2(\mathbb{A}_f)$ we have 
$$ \mathbf{B}(\Q) \cap g U g^{-1} \subset g V g^{-1},$$
However the projection of $\mathbf{B}(\Q) \cap g U g^{-1}$ to the toral $\Q^*$
is a finite subgroup of $\Q^*$, thus contained in $\{\pm 1\}$.  It follows that an index $2$ subgroup of the left-hand side
is contained in $\mathbf{N}(\Q) \cap g U g^{-1}$, which is certainly contained in $g V g^{-1}$
because any open compact of $\mathbf{N}(\Q_q)$ is pro-$q$.  
\qed

\begin{lemma} \label{bfHvanishes}
Suppose that, as above, $V$ is a normal subgroup of $U$.
Let $f=f_{UV}$ be as in (ii) above. 
 Let $\mathcal{F}$ be any sheaf of $\mathcal{O}_{X(V)}$-modules on $X(V)$,
equipped with a compatible action of $U/V$. 

Then
\begin{itemize}
\item[(i)]  For  each $x \in X(V)$,
 the higher cohomology of the stabilizer $(U/V)_x$
on $\mathcal{F}_x$ is  trivial.

\item[(ii)] For each $y \in X(U)$, the higher cohomology of $(U/V)$
acting on $(\pi_* \mathcal{F})_y$ is trivial. 
\end{itemize}

\end{lemma}

\proof 
Note that we can reduce (i) to the case when $(U/V)_x = (U/V)$ by shrinking $U$. 
Both (i) and (ii) will follow, then, if we prove that for any $U/V$-stable affine set $\Spec(A) \subset X(V)$, 
\begin{equation} \label{UVC} \mbox{ higher cohomology of $U/V$  on $ \Gamma(  \Spec(A),  \mathcal{F}) = 0$}. \end{equation}
since the stalks appearing in (i) and (ii) are direct limits of such spaces.

Let $\Delta = U_1/V$ be a Sylow $p$-subgroup of $U/V$; it is sufficient to make the same verification for the higher cohomology of $\Delta$.
Write $B = A^{\Delta}$.  The map
$\Spec(A) \rightarrow \Spec(B)$ is finite {\'e}tale with  Galois group $\Delta$, by (iii) above.   It is now sufficient to show:

\begin{quote}
If $M$ is an $A$-module,  equipped with a $\Delta$-action compatible with its module structure, then $H^q(\Delta, M) =0$ for $q > 0$. 
\end{quote}

Let  $M' = M \otimes_B A$; define a $\Delta$-action  on $M'$ using $g (m \otimes a) = gm \otimes a$
for $g \in \Delta$. 
Since $A$ is a  flat $B$-module, the natural map $H^q(\Delta, M) \otimes_{B} A \rightarrow H^q(\Delta, M')$ is an isomorphism.
We shall show $H^q(\Delta, M') = 0$; the vanishing of $H^q(\Delta, M)$ follows from faithful flatness.

Now 
$M'$ is a  module over $A \otimes_{B} A \simeq \prod_{\delta \in \Delta} A$, and this module structure
is compatible with the $\Delta$-action on $\prod_{\delta \in \Delta} A$ which permutes the factors. 
Therefore, $M'$ is   induced (as a $\Delta$-module) from a representation of the trivial group, and thus has vanishing higher $\Delta$-cohomology by Shapiro's lemma. 
\qed

\subsection{Godement resolution}
Let $T$ be the ``Godement functor'', which assigns to a sheaf $\mathcal{F}$
the sheaf $U \mapsto \prod_{x \in U} \mathcal{F}_x$ of discontinuous sections.  It carries a sheaf of $\mathcal{O}$-modules
to another sheaf of $\mathcal{O}$-modules. 

We will need to discuss the behavior under images. Suppose given a map $f: X' \rightarrow X$ of  schemes.
There is a map of functors
$$ T \rightarrow f_*  T f^{-1}.$$
For a sheaf $\mathcal{F}$ on $X$ and an open set $V \subset X$, this  is given by the natural pullback of discontinuous sections
$$ \prod_{x \in V} \mathcal{F}_x \rightarrow \prod_{x' \in f^{-1} V} (f^{-1} \mathcal{F})_{x'}.$$
If we are working with sheaves of $\mathcal{O}$-modules, then, composing with the natural $f^{-1} \rightarrow f^*$, we get
$T \rightarrow f_* T f^*$, or, what is the same by adjointness,
a natural transformation
$$f^* T \longrightarrow T f^* \mbox{ and (by iterating)}  f^* T^k \rightarrow T^k f^*  .$$
 
 In particular, for a sheaf $\mathcal{F}$ on $X$, there is a map
\begin{equation} \label{rrr} f^* \left(\mbox{ Godement resolution of $\mathcal{F}$}\right) \longrightarrow \mbox{Godement resolution of $f^* \mathcal{F}$.}\end{equation}
 This gives rise to the pullback map in cohomology $H^*(X, \mathcal{F}) \rightarrow H^*(X', f^* \mathcal{F})$.  
 
\subsection{}  \label{MUdef}
 It follows from Lemma \ref{bfHvanishes}
that (with notations as in that Lemma and) for any sheaf $\mathcal{F}$ of $\mathcal{O}_{X(V)}$-modules, 
\begin{equation} \label{UVC3} H^p(U/V, \Gamma(X(V), T \mathcal{F})) = 0, p > 0.\end{equation}
Indeed group cohomology commutes with products (even infinite ones).

 Now let  $\mathcal{G}^{\bullet}(U)$ be the  Godement resolution of $\CL_U$. It is a complex of sheaves
 of $\mathcal{O}_{X(U)}$-modules on $X(U)_S$. 
 Let $M^{\bullet}(U)$ be the global sections of $\mathcal{G}^{\bullet}(U)$: this is a complex of $S$-modules. 
If $V \subset U$, there is a natural action of $U/V$ on $M^{\bullet}(V)$.  It follows from \eqref{UVC3} that

\begin{lemma} \label{acyclicity}
For each degree $i$,  the $U/V$-cohomology of $M^i(V)$ vanishes, i.e. $H^p(U/V, M^i(V))=0$ for $p>0$.
\end{lemma}

 The following result is the crucial one for us. 

\begin{lemma} \label{QI1}
The map arising from \eqref{rrr}
\begin{equation} \label{n m} \mathcal{G}^{\bullet}(U) \rightarrow  \left( f_*  \mathcal{G}^{\bullet}(V)\right)^{U/V}  \end{equation}
 (where $U/V$ denotes invariants)  induces on global sections 
 a quasi-isomorphism
 \begin{equation} \label{tohokuiso} M^{\bullet}(U) \longrightarrow M^{\bullet}(V)^{U/V},\end{equation}
 \end{lemma}
\proof

It is enough to verify that \eqref{n m} is a quasi-isomorphism:
the sheaves $\mathcal{G}^{\bullet}(U)$ and $f_*  \mathcal{G}^{\bullet}(V)^{U/V}$ are flasque --
the latter follows just by examining the definition of the Godement functor $T$ --
and so taking global sections will preserve the quasi-isomorphism.  

Consider the following diagram:
 \begin{equation*}
    \xymatrix{
 \CL_U  \ar[r]^{\sim} \ar[d]^{\sim} & \mathcal{G}^{\bullet}(U) \ar[d]\\ 
  (f_* \CL_V)^{U/V}   \ar[r]^j &  \left( f_*  \mathcal{G}^{\bullet}(V)\right)^{U/V}. 
   }
  \end{equation*}
The left vertical arrow is a quasi-isomorphism:
we have an isomorphism $f_* \CL_V \simeq \CL_U \otimes f_* \mathcal{O}_V$,
and $(f_* \mathcal{O}_V)^{U/V} = \mathcal{O}_U$. 
The  top horizontal arrow is also a quasi-isomorphism. 
It then suffices to show that the arrow $j$ is also a quasi-isomorphism.

The complex $f_* \mathcal{G}^{\bullet}(V)$ is a resolution of $f_* \CL_V $ 
because $f_*$ has no higher cohomology on the quasi-coherent sheaf $\CL_V$. 
Next the stalks of $f_* \CL_V$ and $f_* \mathcal{G}^{\bullet}(V)$ have vanishing $U/V$-cohomology by Lemma \ref{bfHvanishes}. 
Given an acyclic complex of $U/V$-modules supported in degrees $\geq 0$,  each of which have no higher $U/V$-cohomology, 
the $U/V$-invariants remain acyclic. %
This implies that $f_*  \mathcal{G}^{\bullet}(V)^{U/V}$ is a resolution of $(f_* \CL_V)^{U/V}$
as desired.  \qed

\subsection{Compatibility with traces}
We must also mention the compatibility with trace maps. 
Suppose we are given a subgroup $U'$ intermediate between $U$ and $V$:
$$V \subset U' \subset U.$$
We don't require that $U'$ be normal.

There is a natural trace map
$$H^*(X(U'), \CL_{U'}) \rightarrow H^*(X(U), \CL_U).$$

Explicitly the trace $f_* \mathcal{O}_{X(U')} \rightarrow \mathcal \mathcal{O}_{X(U)}$ induces 
$$H^*(X(U'), \omega_{U'}) = H^*(X(U), f_* \omega_{U'})= H^*(X(U), \omega_U \otimes f_* \mathcal{O}_{X(U')}) \stackrel{\mathrm{tr}}{\rightarrow} H^*(X(U), \omega_U).$$ 

 With reference to the identifications of the previous lemma,  this trace map is induced at the level of cohomology by
$$M^{\bullet}(U')  \rightarrow M^{\bullet}(V)^{U'/V} \stackrel{\mathrm{T}}{\rightarrow} M^{\bullet}(V)^{U/V} \stackrel{\sim}{\leftarrow} M^{\bullet}(U)$$
where $\mathrm{T} \in S[U/V]$ is the sum of a set of coset representatives for $U/U'$.

\subsection{Derived invariants and the derived Hecke algebra} 
As in \S \ref{MUdef}, $M^{\bullet}(U)$ is a Godement complex computing the complex of $\omega_U$ 
Now set $$ M_{\infty}^{\bullet}  =  \varinjlim M^{\bullet}(U),$$
which is now a complex of $S$-modules equipped with an action of $G=\GL_2(\Q_q)$. 

We will argue  that the ``derived invariants'' of $U$ on $M_{\infty}^{\bullet}$ gives a complex that computes the cohomology of $X(U)_S$.
We first recall the notion of derived invariants, and its relationship with the derived Hecke algebra.

Let $U$ be an open compact subgroup of $G$. 
Let $U_1 \subset U$ be a normal subgroup with the property that
the pro-order of $U_1$ is relatively prime to $p$. 
 Let $\QQ$ be a projective resolution of $S$ in the category of $S[U/U_1]$-modules;
we regard this as a complex with degree-increasing differential concentrated
in degrees $\leq 0$:
$$\cdots \rightarrow Q_{-2} \rightarrow Q_{-1} \rightarrow Q_0 = S.$$
We may of course regard $\QQ$ as a complex of $S[U]$-modules.

Let 
$\PP = \ind_U^G \QQ$. This is a projective resolution of the smooth $S[G]$ module $S[G/U]$ (in the category of smooth $S[G]$ modules). 
 For any complex $R^{\bullet}$ of $G$-modules, we define the derived $U$-invariants to be the complex
$$\Hom_{S[G]}(\PP, R^{\bullet}) = \Hom_{S[U]}(\QQ, \left( R^{\bullet} \right)^{U_1}).$$
Explicitly, this is a complex whose cohomology computes the hypercohomology $\mathbb{H}^*(U, R^{\bullet})$.

In the case above,  the derived invariants of $U$ on $M^{\bullet}_{\infty}$ compute the cohomology of $X(U)$, in the following sense: 

\begin{lemma} \label{derivedinvariants} 
The natural inclusion of $M^{\bullet}(U) \hookrightarrow M^{\bullet}_{\infty}$
and the augmentation $\QQ \rightarrow S$  induce a quasi-isomorphism:
\begin{equation} \label{brb} M^{\bullet}(U) \stackrel{\sim}{\rightarrow} \Hom_{S[U]}(\QQ, M^{\bullet}_{\infty}) = \Hom_{S[G]}(\PP, M^{\bullet}_{\infty}).\end{equation}
\end{lemma}

\proof

Using  the remarks after \eqref{tohokuiso}, we see that
\begin{equation} \label{above-a} M^{\bullet}(U_1)     \stackrel{\mathrm{q.i.}}{\longrightarrow}   \varinjlim_{U' \subset U_1}  M^{\bullet}(U')^{U_1/U'} \stackrel{\sim}{\rightarrow}  \mbox{$U_1$-invariants on $M_{\infty}^{\bullet}$} .\end{equation}
(for the second arrow: since $U_1$ is  prime to $p$ the functor of taking $U_1$ invariants
commutes with taking a direct limit of smooth $S[U_1]$-modules).

The inclusion $M^{\bullet}(U) \hookrightarrow M^{\bullet}(U_1)$ and the homomorphism $\mathbf{Q} \rightarrow S$  induce
$$ M^{\bullet}(U) \rightarrow \Hom_{U/U_1}(S, M^{\bullet}(U_1)) \rightarrow \Hom_{U/U_1}(\mathbf{Q}, M^{\bullet}(U_1)) $$
and it remains to show that this composite is a quasi-isomorphism.  

The first map is a quasi-isomorphism by Lemma 
\ref{QI1}. To show that the second 
map is a quasi-isomorphism,  it is enough (by a devissage) 
to show that for each fixed degree $j$ 
$$ \Hom(S, M^{j}(U_1)) \rightarrow \Hom(\mathbf{Q}, M^j(U_1))$$
induces a quasi-isomorphism. But the right hand side computes the $U/U_1$ cohomology of $M^j(U_1)$, and we have seen (Lemma \ref{acyclicity}) that this is concentrated in degree 
 zero, where it is just the $U/U_1$ invariants, as needed. 
\qed 

Now we may imitate all the reasoning above, with the role of $\omega$ replaced by $\mathcal{O}$.
Let $N^{\bullet}$ be the corresponding complex. Reasoning as in Lemma \ref{derivedinvariants},  we get a quasi-isomorphism
\begin{equation} \label{Gode} N^{\bullet}(U) \simeq \Hom_{S[U]}(\QQ, N^{\bullet}_{\infty}).\end{equation}
The identification of $S$ with global sections of $\mathcal{O}_{X(V)}$
induce compatible maps $S \rightarrow N^{\bullet}(V)$ for each level $V$, and so 
by passage to the limit a map $S \hookrightarrow N^{\bullet}_{\infty}$.
This induces
\begin{equation} \label{above-b} H^*(U, S) \longrightarrow H^*(X(U), \mathcal{O}).\end{equation}
For $\alpha \in H^j(U, S)$ write $\langle \alpha \rangle \in H^j(X(U), \mathcal{O})$
for its image under this map.  Then we have:

 \begin{lemma} \label{cup product lemma}
Under the identification  $H^*(X(U), \omega_U)$ with the hypercohomology  $\mathbb{H}^*(U, M^{\bullet}_{\infty})$,
(as in the prior Lemma), cup product with $\langle \alpha \rangle$ in Zariski cohomology
is carried to cup product with $\alpha$ in hypercohomology.
\end{lemma} 
\proof

The product $\mathcal{O} \otimes \omega_U \rightarrow \omega_U$
extends to a map $N^{\bullet}(U) \otimes M^{\bullet}(U) \rightarrow M^{\bullet}(U)$ (see \cite[Chapter 6]{Godement}),
which computes on cohomology the cup product. 
This exists compatibly at every level, and  by passage to the direct limit, we arrive at a map 
$ N^{\bullet}_{\infty} \otimes M^{\bullet}_{\infty} \rightarrow M^{\bullet}_{\infty}$ (the tensor product can 
be passed through the direct limit, by \cite[Chapter 2, Prop. 7, \S 6.3]{BourbakiLinearAlgebra}).

Fix a quasi-isomorphism of $S[U]$-modules:
$$q:  \mathbf{Q} \rightarrow \mathbf{Q} \otimes_S \mathbf{Q}.$$

Consider the following diagram, with commutative squares: {\small \begin{equation*}
    \xymatrix{
m' \otimes \alpha'' \in H^i(M^{\bullet}(U)) \otimes  H^j( N^{\bullet}(U)) \ar[d]\ar[rr] && H^{i+j}(M^{\bullet}(U)) \ar[d]  \\ 
m \otimes \alpha' \in  \Hom^i(\mathbf{Q}, M^{\bullet}_{\infty})  \otimes  \Hom^j(\mathbf{Q}, N^{\bullet}_{\infty}) \ar[r]^{\  \qquad \otimes}  &  \Hom^{i+j}(\mathbf{Q} \otimes \mathbf{Q}, M^{\bullet}_{\infty} \otimes N^{\bullet}_{\infty})  \ar[r]^{ \ \ q} &  \Hom^{i+j}(\mathbf{Q}, M^{\bullet}_{\infty}), \\
  m \otimes \alpha \in \Hom^i(\mathbf{Q}, M^{\bullet}_{\infty}) \otimes\Hom^j(\mathbf{Q},  S) \ar[r]^{\ \qquad  \otimes} \ar[u]  & \Hom^{i+j}(\mathbf{Q} \otimes \mathbf{Q}, M^{\bullet}_{\infty} \otimes S)  \ar[r]^{\ \ q}\ar[u] & \Hom^{i+j}(\mathbf{Q}, M^{\bullet}_{\infty}) \ar[u]. 
  }
  \end{equation*}
  }
   where:
  \begin{itemize}
  \item $\Hom$ means in every case homomorphisms of chain complexes of $S[U]$ modules, taken modulo chain homotopy; 
  \item  $\otimes$ comes from the tensor product, which induces a bifunctor on the homotopy category of chain complexes.
  
  \item We fix $m   \in \mathrm{Hom}^i (\mathbf{Q}, M^{\bullet}_{\infty})$, amd 
  $m'$ is the cohomology class corresponding to $m$ under the quasi-isomorphism \eqref{brb}.
  \item We identify $\alpha$ with a class in $\Hom^j(\mathbf{Q},S)$ and  $\alpha'$ is the image of this class, under $S \rightarrow N^{\bullet}_{\infty}$.  Also $\alpha''$ is a cohomology class  in $H^j(N^{\bullet}(U))$ that matches with $\alpha'$ under the quasi-isomorphism \eqref{Gode}.
   \end{itemize}

      The image of $m \otimes  \alpha$, under the bottom horizontal arrows,
  computes the cup product of $m$ and $\alpha$ in $U$-hypercohomology.  This corresponds to the image of  $m \otimes \alpha'$
  in the middle horizontal row.  Finally, this corresponds to the  image of $m' \otimes \alpha''$ in the top row, which gives the Zariski product. \qed

 \subsection{Derived Hecke algebra}
 
Let notation be as above, but specialized to the case $U = K$,
a maximal compact of $\GL_2(\Q_q)$.  
 We may form the differential graded algebra
 $\End_{S[G]}(\PP)$  whose cohomology we understand to be the  (graded) derived Hecke algebra for the pair $(G,K)$. 
   There is an isomorphism (\cite{V3}, (148)) 
 \begin{equation}\label{cosets} \End_{S[G]}(\PP,\PP) ~ \simeq ~ \oplus_{x \in K\backslash G/K} \Hom_{K_x}(\QQ,\QQ_x) \end{equation}
 where $\QQ_x$ is the complex $\QQ$ but with the twisted action of $K_x = K \cap \Ad(g_x) K$ defined by
 $\kappa \ast q = (Ad(g_x^{-1})\kappa)q$; here we have implicitly chosen coset representatives $g_x K$ for each $x \in K \backslash G/K$. 
Taking cohomology, one finds that, for any $i$ there is an isomorphism (\cite{V3}, (149))
 \begin{equation}\label{cosets1}  H^i(\End_{S[G]}(\PP,\PP)) ~~ \isoarrow ~~ \oplus_{x \in K\backslash G/K} H^i(K_x,S) \end{equation}

Now the differential graded algebra $\End_{S[G]}(\PP, \PP)$ acts on 
 $\Hom_{S[G]}(\PP,M_{\infty}^{\bullet})$. Passing to cohomology and applying
 Lemma \ref{derivedinvariants}, we get a graded action of the derived Hecke algebra
for $(G,K)$ 
on $H^*(X_K,\omega_K)$. 
This action is specified by specifying, for each $x  = K g_x  K \in K \backslash G/K$ as above, the corresponding action
of $H^*(K_x, S)$ on coherent cohomology.  
  We can now restate Lemma A.10 of \cite{V3}:
 \begin{lemma}\label{A10}  The action of $h_x \in H^*(K_x,S)$ on $\mathbb{H}^*(K,M_{\infty}^{\bullet})$ is given explicitly by the following composite:
 {\small $$
\xymatrix{ \mathbb{H}^*(K,M_{\infty}^{\bullet}) \ar[r]^{\Ad(g_x^{-1})^*} &  \mathbb{H}^*(K_x,M_{\infty}^{\bullet})  \ar[r]^{m \mapsto g_xm}  & \mathbb{H}^*(K_x,M_{\infty}^{\bullet}) \ar[r]^{m \cup h_x}  &
  \mathbb{H}^*(K_x,M_{\infty}^{\bullet}) \ar[d]^{\mathrm{Cores}} \\
  & & &   \mathbb{H}^*(K,M_{\infty}^{\bullet}) 
}
  $$}
  \end{lemma}

We obtain the derived Hecke operator $T_{q,z}$ described in \S \ref{DHmain},
with the coefficient ring $S=\mathcal{O}/\mathfrak{p}^n$, 
by taking $x = \left(\begin{array}{cc} q & 0 \\ 0 & 1 \end{array}\right)$
and by taking the cohomology class $h_x \in H^1(K_x, S)$  as the composite:
$$ \left( \begin{array}{cc} a & b \\ c  & d \end{array}\right) \in K_x \mapsto \langle a/d \mbox{ mod } q, z \rangle,$$
where $z \in k \langle -1 \rangle$ is regarded as a homomorphism $(\Z/q)^* \rightarrow \mathcal{O}/\mathfrak{p}^n$. 
Indeed, to verify this, it only remains to show that  the induced map  \begin{equation} \mathbb{H}^*(K_x,M_{\infty}^{\bullet})  \rightarrow \mathbb{H}^*(K_x, M_{\infty}^{\bullet})\end{equation}
given by cupping with the class $h_x $ is identified with 
$$ H^*(X(K_x), \omega) \stackrel{\cup z \Shim}{\longrightarrow} H^*(X(K_x), \omega),$$ 
that is to say the cup 
product with  $z\Shim$, i.e. the Shimura class multiplied by $z$, regarding
as a class in the cohomology of $X(K_x)$ with coefficients in $\mathcal{O}/\mathfrak{p}^n$. 
This follows easily from Lemma \ref{cup product lemma}.

The following remark is due entirely to the first-named author (M.H); the second-named author disclaims both credit and responsibility for it.

 \begin{rmk}\label{millenials}  For the benefit of those millenials who believe the Godement resolution is one of the founding documents of the United Nations, here is a translation of the above construction into contemporary language.  We thank Nick Rozenblyum for his patient guidance.  We work in the DG category (or stable $\infty$-category) $C$ of complexes of quasicoherent sheaves on the scheme $X_{\infty}$, and consider the object $\CL_{\infty}$, all over $Spec(S)$.  This object carries an action by $G = GL_2(\QQ_q)$.  Therefore the object $R\Gamma(\CL_{\infty})$ in the DG category $Mod_S$ of complexes of $S$-modules carries an action of $G$.  Everything up through Lemma \ref{cup product lemma} is automatic in this setting.  The remaining observations are not strictly necessary to formulate the conjecture; however, they do provide the explicit computation of the derived Hecke operator, as in Lemma \ref{A10}, needed in order to test the conjecture in specific applications.
\end{rmk}

\section{Magma Code}

What follows is a sample of Magma code which we used to compute the derived Hecke operator for the modular form of level $31$,
with $q=139$ and $p=23$. 

\begin{alltt}
N := 31;
Q := 139;
L := 23;
F := FiniteField(L);
M := ModularForms(N*Q);
S := CuspidalSubspace(M);
SQ := BaseExtend(S, RationalField()); 
SF := BaseExtend(S, F); 
V, h := VectorSpace(SF); 
time Tq := HeckeOperator(SF,N);
time Wq := AtkinLehnerOperator(SF,N);
Iq := IdentityMatrix(F, Dimension(S));
Qq := Iq +Wq*Tq; /* Qq projects from level QN back down to level Q */
Pro := Dimension(S);
Z<q> := PowerSeriesRing(IntegerRing());
QQ<q> := PowerSeriesRing(RationalField());

CUTOFF := Dimension(S)+3; 

eps := KroneckerCharacter(-N); 
WeightOneSpace := ModularForms(eps, 1); 
etatemp := WeightOneSpace.2;
etaprodA := qExpansion(etatemp, CUTOFF);
etaprodB := Composition(etaprodA, q^Q+O(q^CUTOFF));
g := etaprodA * etaprodB + O(q^CUTOFF);
 
g0 :=  SF ! g;  
 W := Vector(F, Inverse(h)(g0));
Wfin := W * Qq;
/*denom := Denominator(Wfin);
print(Factorization(denom)); */
M2 := ModularForms(Q);
S2 := CuspidalSubspace(M2);
S2Q := BaseExtend(S2, RationalField());
S2F := BaseExtend(S2, F); 
V2,h2 := VectorSpace(S2F); 
CUTOFF2  := Dimension(S2)+3;
 projform := S2F  ! h(Wfin);  
projformcoeff := Vector(F, Inverse(h2)(projform)); 
 normcoeffF := projformcoeff;  
randprime := 41;
randT := HeckeOperator(S2F, randprime);      
charpoly := CharacteristicPolynomial(randT); 
P<u>,h3 := ChangeRing(PolynomialRing(IntegerRing()), F); 
 Factorization(P! charpoly); 
unnormalizedredpoly := charpoly/(u-randprime-1); 
redpoly := unnormalizedredpoly/Evaluate(unnormalizedredpoly,randprime+1);
print(Evaluate(redpoly, randprime+1));
projmatrix := Evaluate(P! redpoly, randT);
 finalanswerinbasis := normcoeffF *  projmatrix;
print(finalanswerinbasis); 
\end{alltt}


\begin{thebibliography}{99}
\bibitem{SGA4}  M. Artin, A. Grothendieck, and J. Verdier, eds. S\'eminaire de G\'eom\'etrie Alg\'ebrique du Bois Marie - 1963-64 - Th\'eorie des topos et cohomologie \'etale des sch\'emas - (SGA 4) - vol. 3, {\it Lecture Notes in Mathematics}, {\bf 305},  Berlin; New York: Springer-Verlag. (1972). 

\bibitem{ArtinMilne}  M. Artin, J. S. Milne, Duality in the flat cohomology of curves, {\it Invent. Math.}, {\bf 35} (1976) 11--129. 

 \bibitem{Calegari} F. Calegari, Non-minimal modularity lifting in weight one, {\it J. Reine angew. Math.}, Published online  2015-12-17.

\bibitem{CG} F. Calegari and D. Geraghty, Modularity lifting beyond the Taylor-Wiles method. Preprint available on arxiv.

\bibitem{DR}  P. Deligne and M. Rapoport, Les sch\'emas de modules de courbes elliptiques, in {\it Modular functions of one variable II.}  Springer Berlin Heidelberg, (1973) 143--316.

\bibitem{DS}   P. Deligne and J.-P. Serre, .Formes modulaires de poids $1$, {\it Annales scientifiques de l'\'Ecole Normale Sup\'erieure}, {\bf 7},  (1974). 507--530.


\bibitem{GV}  S. Galatius and A. Venkatesh, Derived Galois deformation rings. Preprint available on arxiv.

\bibitem{Godement} R. Godement, {\it Topologie alg\'ebrique et th\'eorie des faisceaux}, Paris:  Hermann (1958).

\bibitem{tohoku}  A. Grothendieck, Sur quelques points d'alg\`ebre homologique, {\it Tohoku Mathematical Journal}  Second Series, {\bf 9} (1957) 119--183.


\bibitem{groth}  A. Grothendieck; Le groupe de Brauer. III. Exemples et compl\'ements. {\it Dix Expos\'es sur la
Cohomologie des Sch\'emas}, North-Holland, Amsterdam (1968) 88--188.

\bibitem{H}  M. Harris, The Taylor--Wiles method for coherent cohomology, {\it J. reine  angew. Mathematik}, {\bf 679}  (2013) 125--53.

\bibitem{K}  N. Katz, p-adic properties of modular schemes and modular forms, {\it Modular functions of one variable, III }(Proc. Internat. Summer School, Univ. Antwerp, Antwerp, 1972), {\it Lecture Notes in Mathematics}, {\bf 350}, Berlin, New York: Springer-Verlag, (1973) 69--190. 

\bibitem{mazur} B. Mazur, Modular curves and the Eisenstein ideal. {\it Inst. Hautes \'Etudes Sci. Publ. Math.}
{\bf 47} (1977) 33-- 186.

\bibitem{Merel}  L. Merel, L'accouplement de Weil entre le sous-groupe de Shimura et le sous-groupe cuspidal de $J_0(p)$, {\it J. Reine Angew. Math. } {\bf 477} (1996), 71--115. 

\bibitem{PV}  K. Prasanna and A. Venkatesh, Automorphic  cohomology, motivic cohomology, and the adjoint $L$-function. Preprint available on arxiv.


\bibitem{Scholl}  A. Scholl, Integral elements in K-theory and products of modular curves
In: {\it The Arithmetic and Geometry of Algebraic Cycles} ed. B. B. Gordon et al. NATO Science Series C, volume 548 (Kluwer, 2000), 467--489.

\bibitem{V3}  A. Venkatesh, Derived Hecke algebra and cohomology of arithmetic groups. Preprint available on arxiv.


\end{thebibliography}
 \end{document}